\documentclass[11pt]{book}

\usepackage{amsmath}
\usepackage{latexsym}
\usepackage{amsfonts}
\usepackage{graphics}
\usepackage{amssymb}

\define\R{\Bbb{R}}

\define\Z{\Bbb{Z}}

\begin{document}



\newtheorem{theorem}{Theorem}[chapter]
\newtheorem{lemma}{Lemma}[chapter]
\newtheorem{definition}{Definition}[chapter]
\newtheorem{claim}{Claim}[chapter]
\newtheorem{corollary}{Corollary}[chapter]
\newtheorem{condition}{Condition}[chapter]
\newtheorem{question}{Question}[chapter]
\newtheorem{example}{Example}[chapter]
\newtheorem{remark}{Remark}[chapter]

\newbox\qedbox
\setbox\qedbox=\hbox{q.e.d.}

\newenvironment{proof}{\smallskip\noindent{\bf Proof.}\hskip \labelsep}
                        {\hfill\penalty10000\copy\qedbox\par\medskip}

\mainmatter{

\begin{center}{\bf CLASSIFICATION OF STABLE MINIMAL SURFACES 
                 BOUNDED BY JORDAN CURVES IN CLOSE PLANES}
\end{center}
\vskip .5cm

 \centerline{ROSANNA PEARLSTEIN}
\vskip 1cm

\centerline{\bf Abstract} 

      We study compact stable embedded minimal surfaces 
      whose boundary is given by two collections of closed
      smooth Jordan curves in close planes of Euclidean 
      Three-Space.  Our main result is a classification 
      of these minimal surfaces, under certain natural 
      geometric asymptotic constraints, in terms of certain 
      associated varifolds which can be enumerated explicitely.
      One consequence of this result is the uniqueness of the 
      area-minimizing examples.  Another is the asymptotic 
      non-existence of stable compact embedded minimal surfaces 
      of positive genus bounded by two convex curves in parallel 
      planes.

\vskip 1cm

\centerline{\bf Acknowledgments}

The author wishes to thank her thesis advisor, professor 
Bill Meeks, for suggesting this problem to her, and for 
many helpful conversations.
 
\vskip 1cm

\setcounter{chapter}{1}
\setcounter{section}{-1}
\section{Introduction}              

At the Princeton Bicentennial Conference in $1946$, Tibor Rad{\'o}
posed the problem to estimate the number of compact minimal surfaces
bounded by one or more given Jordan curves in Euclidean three-space.
We will answer this question of Rad{\'o} under the assumption that the
surfaces are embedded and stable, and their boundary curves satisfy a
natural asymptotic geometric constraint.  We now describe our results
and methods in more detail.  Throughout this paper, $\cal A$ and $\cal
B$ will each denote a collection of disjoint closed smooth Jordan
curves in the plane ${P_0}:z=0$ of ${\R}^3$.  We make the assumption
that $\cal A$ and $\cal B$ intersect transversely in a finite
collection of points, which we call {\em crossing points}.  We denote
by ${\cal B}(t)$ the vertical translation of $\cal B$ to the plane
${P_t}:z=t$.  Consider the immersed $1$-cycle $Z={\cal A}{\cup}{\cal
B}$ in $P_0$, and let ${\cal C}({\cal A}, {\cal B})$ denote the
collection of the closures of the bounded components in
${P_0}{\backslash}Z$.  Let ${\cal V}({\cal A}, {\cal B})$ denote the
(finite) collection of integer multiplicity varifolds that can be
represented by nonnegative integer multiple sums of the components in
${\cal C}({\cal A}, {\cal B})$, where the multiplicity changes by
$one$ in passing from any component to an adjacent one, and at least
one of the four components meeting at a crossing point has
multiplicity zero (note that as a consequence these varifolds have $Z$
as their ${\Z}_2$-boundary).  Let ${\cal S}(t)$ be the collection of
compact stable embedded minimal surfaces bounded by ${\cal
A}{\cup}{\cal B}(t)$.  Our main theorem is:

\begin{theorem}[Main Theorem]  For values of $t$ sufficiently    small,
one can  naturally associate to a  varifold $V$ in ${\cal V}({\cal A},
{\cal B})$  a  unique   compact   stable  embedded  minimal    surface
${\Sigma}(V, t)$ in   ${\cal  S}(t)$ bounded by ${\cal   A}{\cup}{\cal
B}(t)$.   Furthermore, this association  between  varifolds in  ${\cal
V}({\cal A}, {\cal B})$ and surfaces in  ${\cal S}(t)$ is a one-to-one
correspondence. \end{theorem}

In order to discuss further the geometry and topology of the stable
minimal surface ${\Sigma}(V, t)$ associated to a varifold $V{\in}{\cal
V}({\cal A}, {\cal B})$, we need to give a few definitions.  Note that
around each crossing point, $Z$ divides a sufficiently small
neighborhood $D$ in four connected components; we order them according
to the counterclockwise orientation of $D$, so that the multiplicity
of the first component is zero.  Now fix a $V{\in}{\cal V}({\cal A},
{\cal B})$.  Our assumptions on $V$ imply that only the multiplicities
$(0, 1, 0, 1)$ and $(0, 1, 2, 1)$ can occur at the crossing points.
We let $v_1$ be the number of crossing points with multiplicity $(0,
1, 0, 1)$ in ${\cal A}{\cap}{\cal B}$, let $v_2$ be the number of
crossing points with multiplicity $(0, 1, 2, 1)$ in ${\cal
A}{\cap}{\cal B}$, let $f_1$ be the number of components of $V$ with
multiplicity one, let $f_2$ be the number of components of $V$ with
multiplicity two, let $e_1$ be the number of multiplicity one
connected components of $Z{\backslash}{\{}crossing\ \ points{\}}$, let
$e_2$ be the number of multiplicity two connected components of
$Z{\backslash}{\{}crossing\ \ points{\}}$.  Then we have the following
result.

\begin{theorem} Let $V{\in}{\cal V}({\cal A}, {\cal B})$, and for $t$
sufficiently small let ${\Sigma}(V, t)$ be the surface given in the
statement of theorem $1.1$.  Then the Euler characteristic of
${\Sigma}(V, t)$ is equal to $({v_1} + 2{v_2}) - ({e_1} + 2{e_2}) +
({f_1} + 2{f_2})$.  \end{theorem}

It follows  that since there is  exactly one varifold ${V_0}{\in}{\cal
V}({\cal A}, {\cal B})$ having least  area (namely no multiplicity $2$
components),   there is exactly   one  ${\Sigma}({V_0}, t)$ in  ${\cal
S}(t)$ with Euler  characteristic equal to $N -  l + f$,  where $N$ is
the number of crossing points, $l$ the  number of connected components
of  $Z{\backslash}{{}     crossing\   points{}}$, $f$ is the number of
multiplicity one components of $V$.      This     surface
${\Sigma}({V_0}, t)$  is the  unique surface  of  least area in ${\cal
S}(t)$.

The first step in our proof of theorem $1.1$ is to  show that, for any
fixed  varifold $V$, there  exists  a compact  stable embedded minimal
surface bounded by ${\cal  A}{\cup}{\cal B}(t)$ corresponding  to $V$,
when  $t$ is   sufficiently  small.   More   precisely, we  prove  the
following theorem, which follows from   an existence result proven  by
W. Meeks and S. T.    Yau \cite{MeYa1}.  \begin{theorem} Let $\cal  A$
and ${\cal B}(t)$ be two collections of  disjoint closed smooth Jordan
curves in the  planes ${P_0}$ and  $P_t$ respectively, such that $\cal
A$ intersects ${\cal B}={\cal B}(0)$  transversely in a finite  number
of points.   Fix a varifold   $V$ in ${\cal  V}({\cal A},  {\cal B})$.
Then, for every sufficiently small  $t$, there exists a compact stable
embedded minimal surface ${\Sigma}(t)$ with ${\partial}({\Sigma}(t)) =
{\cal A}{\cup}{\cal B}(t)$,  and  such that  ${\Sigma}(t)$ is  ``close
to''  $V$,   in   the sense that    the  family  of   minimal surfaces
${\{}{\Sigma}(t){\}}_t$  converges    to    $V$  as   $t$   approaches
$0$. \end{theorem} Next we describe  geometrically the stable  minimal
surfaces obtained in the previous theorem, extending a previous result
by W.  Rossman \cite{Ro}.

\begin{theorem} Let ${\cal A}{\cup}{\cal B}(t)$ be the boundary of
compact stable embedded minimal surfaces ${\Sigma}(t)$ arising in
theorem $1.3$.  Then, for $t$ sufficiently close to $0$, the surfaces
${\Sigma}(t)$ have the following properties:

\begin{description} 

\item{1.} In the complement of small vertical cylindrical
neighborhoods $N(i)$ of the crossing points, the components of
${\Sigma}(t){\backslash}{\bigcup}N(i)$ are graphs over their
projections to the plane $P_0$.

\item{2.} In $N(i)$, ${\Sigma}(t)$ is either approximately helicoidal,
or is the union of two graphs over the plane $P_0$.

\end{description} 
\end{theorem}

Next we prove that if a minimal surface ${\Sigma}(t)$ possesses a
description as in theorem $1.4$, then such minimal surface must be
stable, for $t$ sufficiently small.

\begin{theorem} Let ${\Sigma}(t)$ be a minimal surface with boundary
${\cal A}{\cup}{\cal B}(t)$.  Suppose that ${\Sigma}(t)$ is described
as in theorem $1.4$.  Then, for $t$ sufficiently close to $0$,
$\Sigma(t)$ is stable. \end{theorem}

The proof of the above theorem involves some technical lemmas; we give
here a sketch of the proof of theorem $1.5$, referring the reader to
section $4$ for details.  A key step of the proof is the following
observation:  let $\Sigma$ be a minimal
surface, given by the union ${{\Sigma}_1} {\cup} S {\cup}
{{\Sigma}_2}$, and suppose that the first eigenvalue of the Jacobi
operator on ${\Sigma}_1$ and ${\Sigma}_2$ is strictly positive, and
$S$ is a simply connected sufficiently flat piece, in the sense that the diameter of
its Gaussian image is very small.  Then the first eigenvalue of the
Jacobi operator on the minimal surface $\Sigma$ is also strictly
positive.   Our hypotheses allow
us to apply some techniques developed by Kapouleas in \cite{Ka}, and
used in his construction of constant mean curvature surfaces in
${\R}^3$. These results imply that if one glues two stable pieces
$E_1$ and $E_2$ by a sufficiently flat bridge $S$, then no Jacobi
vector fields arise on ${E_1} {\cup} S {\cup} {E_2}$, since none arise
from $E_1$ or $E_2$, and no new ones arise from gluing $S$.  In the
case corresponding to a surface ${\Sigma}(V, t)$ in theorem $1.5$, we
start with $E_1$ and $E_2$ each being a small helicoidal minimal
surfaces near a crossing point, and $S$ being one of the flat
components of ${\Sigma}(t){\backslash}{ \bigcup}N(i)$ adjacent to
$E_1$, and $E_2$, and is the region that glues $E_1$ with $E_2$.  We
can do this by the description of ${\Sigma}(V, t)$ given in theorem
$1.4$.  A proof by induction on the number of helicoidal components
shows that ${\Sigma}(V, t)$ is stable.
This completes our brief sketch of theorem $1.5$.

\begin{remark}  This inspires a {\em Gauss-map bridge principle},
which will be discussed more extensively elsewhere.\end {remark}

We now give a brief sketch of the proof of our Main Theorem.  By
theorem $1.3$, we can associate to each varifold $V{\in}{\cal V}({\cal
A}, {\cal B})$ a minimal surface ${\Sigma}(V, t)$, for $t$
sufficiently small.  The Main Theorem states that ${\Sigma}(V, t)$ is
eventually unique.  We now give a sketch of the proof of the
uniqueness.  Suppose that there are two sequences of stable minimal
surfaces ${\Sigma_1}(t_i)$ and ${\Sigma_2}(t_i)$ which have
${\partial}({\Sigma_1}(t_i))={\partial}({\Sigma_2}(t_i))={\cal
A}{\cup}{\cal B}(t_i)$, where ${t_i}$ approaches $0$ as $i$ tends to
infinity.  We show that if such a sequence of minimal surfaces
existed, then there would be two other sequences of embedded stable
minimal surfaces with boundary given by ${\cal A}{\cup}{\cal B}(t_i)$,
and these surfaces could be taken so that their interiors would be
disjoint.  So we may assume that ${\Sigma_1}(t_i)$ and
${\Sigma_2}(t_i)$ are disjoint in their interior, for all $t_i$, and
that both sequences converge to a fixed varifold $V$, as ${t_i} \to
0$.  Theorem $1.4$ implies that, for $i$ large, the surfaces
${\Sigma_1}(t_i)$ and ${\Sigma_2}(t_i)$ bound a product region $R(i)$,
where the interior angles are less than $\pi$, and the surfaces are
ambient isotopic in $R(i)$.  By a deep minimax theorem by Pitts and
Rubinstein \cite{PiRu}, generalized to the case of nonempty boundary
by Jost \cite{Jo}, there would exist an unstable embedded minimal
surface $\Sigma^*(t_i)$ in $R(i)$, such that
${\partial}(\Sigma^*(t_i))={\partial}({\Sigma_1}(t_i))={\partial}({\Sigma_2}(t_i))={\cal
A}{\cup}{\cal B}(t_i)$.  Because ${\Sigma_1}(t_i)$ and
${\Sigma_2}(t_i)$ are expressed as normal graphs over each other, and
by the definition of minimax, we are able to show that, away from the
crossing points, $\Sigma^*(t_i)$ is a very flat graph.  Then Nitsche's
$4{\pi}$-theorem allows us to show that $\Sigma^*(t_i)$ is
approximately helicoidal around the crossing points.  Hence the
unstable minimal surface $\Sigma^*(t_i)$ would possess the same
geometric description as ${\Sigma_1}(t_i)$ and ${\Sigma_2}(t_i)$. Then
theorem $1.5$ implies that, for $i$ sufficiently large,
$\Sigma^*(t_i)$ would have to be stable.  This produces a
contradiction, and finishes the sketch of our proof of the Main
Theorem.

  We then give an upper bound on the number of compact embedded 
 stable minimal surfaces bounded by two given collections of curves.
 More precisely, we show the following.

\begin{corollary} Once the limiting cycle $Z={\cal A}{\cup}{\cal B}$
in the plane $P_0$ is given, the number of stable compact minimal
surfaces ${\Sigma}(t)$ such that ${\partial}({\Sigma}(t)){\to} Z$ as
$t{\to}0$, is bounded above by

$$ 2^{f^i_-} + 2^{f^o_-},
$$for $t$ sufficiently close to $1$, where $f^i_-$ is the number of
bounded components of $P_0{\backslash}Z$ inside
${Region}({\cal A}) {\cap} {Region}{({\cal B})}$, and $f^o_-$
is the number of bounded components of $P_0{\backslash}Z$ outside
${Region}({\cal A}) {\cup} {Region}({\cal B})$.
\end{corollary}

  Finally, we observe how the above results support a conjecture made
by W.H. Meeks in \cite{Me1}, about the nonexistence of positive genus
compact embedded minimal surfaces bounded by two convex curves in
parallel planes of ${\R}^3$.

\setcounter{chapter}{2}
\setcounter{section}{-1}
\setcounter{theorem}{0}
\setcounter{lemma}{0}
\setcounter{definition}{0}
\setcounter{claim}{0}
\setcounter{corollary}{0}
\setcounter{condition}{0}
\setcounter{question}{0}
\setcounter{example}{0}
\setcounter{remark}{0}

\section{Preliminary notions}

In this section we recall some known facts about minimal surfaces.

\section{Minimal surfaces are locally area-minimizing}

There are several possible ways to define minimal surfaces in
${\R}^3$.  A surface $S$ in ${\R}^3$ is {\it minimal}\ \ if for every
point $p$ of $S$ it is possible to find a neighborhood $U_p$ of $p$,
contained in $S$, which is the unique surface of least area among all
surfaces with boundary ${\partial}U_p$.  Notice that $S$ could have
infinite area, for instance, if $S$ is a plane in ${\R}^3$.  A minimal
surface can also be defined as a surface for which every compact
subdomain $C$ is a critical point for the area functional among all
surfaces having boundary equal to ${\partial}C$.  From the point of
view of local geometry, this is equivalent to the condition that the
mean curvature $H$ be identically zero. 

\begin{theorem}[Meusnier] A regular surface $\Sigma$ immersed in 
${\R}^3$ is a critical point for the area functional if and only if 
the mean curvature vanishes identically on $\Sigma$.  \end{theorem}

\section{Stability: eigenvalues of the Laplacian}

Recall that the Gauss map of an oriented surface in ${\R}^3$ assigns
to each point of the surface the unit normal vector to the surface at
that point, viewed as an element of the unit sphere $S^2$ in ${\R}^3$.
We fix an orientation on $S^2$ by the inward pointing normal vector.

\begin{theorem}[Christoffel] Given a surface 
${\Sigma}\hookrightarrow{\R}^3$, the Gauss map $g$ is conformal or 
anti-conformal if and only if $\Sigma$ is a minimal surface 
(or a sphere).  
\end{theorem}

\begin{proof}
  The proof of this theorem follows from the definitions
  of mean curvature and Gauss map.  A complete proof can be
  found in \cite[p. 385, vol. 4]{Sp}, and the original 
  proof in\cite{Ch}.
\end{proof}

Recall that a minimal surface is stable if, for each subdomain $D$ the
second derivative of area is positive for each smooth normal variation
of $D$ which is the identity on ${\partial D} = \Gamma$, and 
{\it unstable} if for some smooth normal variation of $D$ which is the
identity on ${\partial}D = \Gamma$ the second derivative of area is
negative.

We now recall another characterization of stability. We start
with the second variational formula (see for example \cite{Si}):

\begin{theorem}[Second variational formula]
\protect\begin{equation} {\frac{d^2A}{dt^2}}{\bigg|}_{t = 0} =
  {\int_{D}(|\nabla h|^2 + 2Kh^2){\,}dA}, \protect\end{equation}where
$h$ is as in the previous paragraph.
\end{theorem}

Clearly the stability of $D$ is equivalent to the inequality:

$$
\int_{D}|\nabla h|^2dA > -2\int_{D}Kh^2dA
$$for any $h$ in $C^\infty(\overline D)$ which vanishes on
${\Gamma}$.  Recall now that for a minimal surface, the pullback
metric of the metric on $S^2$ under the Gauss map is given by:

\protect\begin{equation} d{\hat s} = -Kds^2,
\protect\end{equation}from which it follows that

$$ d{\hat A} = -KdA
$$and

$$ |\nabla h|^2 = -K|{\hat {\nabla}}h|^2,
$$where ${\hat {\nabla}}$ is the gradient of the new (pullback) metric
on $\Sigma$.  Hence we can rewrite the previous equation in the form:

$$
{\int}_{D}|{\hat {\nabla}}h|^2{\,}d{\hat A}>2{\int}_{D}h^2{\,}d{\hat A}.
$$Now recall that the ratio

$$
Q_{D}(h) = \frac{{\int}_{D}|{\nabla}h|^2{\,}dA}{{\int}_{D}h^2{\,}dA}
$$is the {\it Rayleigh quotient}, and the quantity

$$
{\lambda}_1(D) = \inf Q_{D}(h)
$$represents the first (smallest) eigenvalue of the Dirichlet problem

\[ \left\{ \begin{array}{lll}
          {\Delta}h + {\lambda}h  &= 0  & \mbox{in D} \\
             h  &= 0  & \mbox{on $\Gamma$.} 
             \end{array}
           \right. \]The above infimum can be taken over all
           continuous piecewise smooth functions on the closure of $D$
           which vanish on $\Gamma$, and $\Delta$ represents the
           Laplace operator with respect to the (given) metric on $D$.
           If ${\Gamma}$ is sufficiently smooth, then the above
           boundary value problem has a solution $h_1$ corresponding
           to the eigenvalue ${\lambda}_1$, and the above infimum is
           attained by $h = h_1$.  Putting these facts together, we
           see that the above conditions of stability are equivalent
           to the inequality:

$$ {\lambda}_1(D) > 2.
$$For a minimal surface, the metric $(2.3)$ is exactly the pullback
metric under the Gauss map of the metric on $S^2$.  Thus we have the
following theorem, due to Schwarz:

\begin{theorem}[Schwarz]  Let $D$ be a (relatively) compact domain on 
  a minimal surface ${\Sigma}{\hookrightarrow}{{\R}^3}$.  Suppose that
  the Gauss map $g$ takes $D$ injectively on a domain $\hat{D}$ on the
  unit sphere $S^2$.  If ${\lambda_1}(\hat{D}) < 2$, then $D$ cannot
  be locally area-minimizing with respect to its boundary.
\end{theorem}

\begin{proof}  Because of the above observations we see that the
  condition ${\lambda_1}(\hat{D}) < 2$ implies the existence of a
  $C^{\infty}$ function $u:{\overline{D}}{\to}\R$ such that
  $u_{|{\partial}\hat{D}} = 0$ and which makes the second variation of
  area (2.2) negative.  Namely $D$ is the initial point of a
  1-parameter family of minimal surfaces with boundary equal to
  $\partial{D}$, such that each element in this family has area
  smaller than the area of $D$.  Hence $D$ cannot be locally
  area-minimizing (and therefore cannot be stable).
\end{proof}

Schwarz's theorem, and the facts that the first eigenvalue
of the Laplacian on a hemisphere of $S^2$ is equal to $2$, and that
$D'{\subseteq}D{\Rightarrow}{\lambda_1}(D'){\geq}\lambda_1(D)$,
imply that if the Gauss map $g$ is injective on a domain $D$ of a
minimal surface, and the Gaussian image $g(D)$ contains a hemisphere
on $S^2$, then $D$ is unstable.  Barbosa and Do Carmo
\cite{BaDo1,BaDo2} generalized this result, showing that $g$ does not
necessarily have to be injective, by proving the following theorem.

\begin{theorem}[Barbosa-Do Carmo] If the area of the Gaussian image 
  $g(D){\subseteq}{S^2}$ is less than $2\pi$, then $D$ is stable.
\end{theorem}
 
\begin{remark} Note that the converse statement to theorem $2.5$ 
  does not hold.  This can be seen by considering two or more disjoint
  stable minimal annuli, and by gluing them together by thin bridges,
  to obtain a new connected stable minimal surface (see \cite{MeYa1},
  \cite{Sm}, \cite{Wh2}, \cite{Wh3}).  This new minimal surface is
  stable, but in general its image under the Gauss map, counted with
  multiplicities, has area larger that $2\pi$.
\end{remark}

We recall yet another characterization \cite{BaDo1,BaDo2} of stable
minimal surfaces.  First we need the following definition.

\begin{definition} A smooth normal vector field $V(p) = p + u(p){\vec N}(p)$,
  where ${\vec N}(p)$ is the normal vector to $D$ in $p$ is said to be
  a Jacobi vector field if $u:{\overline D}{\to}{\R}$,
  $p{\mapsto}u(p)$ satisfies the Jacobi equation
  $$-{\Delta}u + 2uK = 0,$$ where $\Delta$ is the Laplace operator on
  $D$ and $K$ is the Gaussian curvature on $D$.
\end{definition}

Then it can be shown, via the Morse Index Theorem \cite{SSm}, that

\begin{theorem}[Classical]
  A domain $D$ is stable if and only if no subdomain $D'$ inside $D$
  admits a Jacobi vector field which is nonzero in $D'$ but zero on
  $\partial {D'}$.
\end{theorem}

\section{The maximum principle}

Suppose that $X=(X_1,X_2,X_3):S{\to}{{\R}^3}$ is a minimal immersion.
Then one can check, using the definition of zero mean curvature
(since ${\Delta} X= 2 H{\cdot}N$), that the coordinate functions $X_1$,
$X_2$, $X_3$ are harmonic functions.  It is well known that there is a
maximum principle for harmonic functions.  Hence there is also a maximum
principle for minimal surfaces.

\begin{theorem}[Maximum principle at an interior point]
  Let $M_1$ and $M_2$ be two minimal surfaces in ${\R}^3$ that
  intersect at an interior point $p$.  If $M_1$ lies on one side of
  $M_2$ near $p$, then $M_1$ and $M_2$ coincide in a neighborhood of
  the point $p$.
\end{theorem}

\begin{theorem}[Maximum principle at a boundary point]  Suppose \newline
  that the minimal surfaces $M_1$ and $M_2$ have boundary curves $C_1$
  and $C_2$ respectively.  Let $p$ be a point in ${C_1}{\cap}{C_2}$,
  and suppose that $T_p(M_1) = T_p(M_2)$ and $T_p(C_1) = T_p(C_2)$.
  Choose orientations of $M_1$ and $M_2$ in such a way that the two
  surfaces have the same normal vector $\vec n$ at $p$.  If near $p$
  $M_1$ lies on one side of $M_2$ and the conormal vectors of $M_1$
  and $M_2$ at $p$ coincide, then $M_1$ and $M_2$ coincide in a
  neighborhood of the point~$p$.
\end{theorem}

\section{Minimax theorems}

It is intuitively evident that a smooth function $f$ from ${\R}^N$ to
$\R$ which has two isolated nondegenerate relative minimum points 
should also have at least one unstable critical point.  This was 
proven by Courant in his book \cite{Co}.

Later in this paper we will need the existence of a {\em minimax
  solution} between two stable minimal surfaces with the same
boundaries (see the Introduction and section $4$).

\begin{definition}  Let $M$ be the space of embedded surfaces contained
  in the product region of space bounded by
  ${\Sigma_1}{\cup}\Sigma_2$, and having interior angles less than
  $\pi$.  Let $\overline M$ be an opportunely constructed
  compactification of $M$ (see \cite[p. 234-235]{Jo}), and let $\cal
  F$ be a cover of $\overline M$ by connected sets.  Denote with $F$ a
  generic element in $\cal F$.  Define now the minimax of the function
  $A =$ ``${Area}$ on $\cal F$'' by:

$${Minimax}(A, {\cal F})={inf}_{F{\in}{\cal F}}{sup}{\{ A(M) | M{\in}F{\}}}.$$

\end{definition}

Then the deep minimax theorems proven by Pitts and Rubinstein
\cite{PiRu}, and generalized by Jost \cite{Jo} to the case of nonempty
boundary, guarantee that the following assertion holds true.

\begin{theorem}[Minimax] There exists an embedded unstable minimal
surface ${\Sigma}^*$, contained in the region of space $\cal R$
bounded by ${\Sigma_1}{\cup}\Sigma_2$, having
${\partial}({\Sigma}^*)={\partial}({\Sigma_1})={\partial}(\Sigma_2)$.
Furthermore, the genus of ${\Sigma}^*$ is at most equal to the genus
of ${\Sigma}_i$, $i = 1, 2$.
\end{theorem}

\begin{remark} In our case, it will be very important to make sure that 
  ${\Sigma}^*$ is homeomorphic to each of ${\Sigma}_1$ and
  ${\Sigma}_2$.  More precisely, the minimax theorems mentioned above
  guarantee that the genus of ${\Sigma}^*$ is at most equal to $g$,
  but it could, {\em a priori}, be less than $g$.  However, since
  ${\Sigma}^*$ is contained in the product region $\cal R$ bounded by
  ${\Sigma_1}{\cup}\Sigma_2$, and both ${\Sigma}_1$ and ${\Sigma}_2$
  are incompressible in $\cal R$, then the genus of ${\Sigma}^*$ is at
  least $g$.  Hence the unstable minimax solution is homeomorphic to
  the two stable ones.
\end{remark}

\section{Nitsche's $4\pi$-theorem}

Finally, let us state here another result which will be useful in
subsequent sections of this paper, namely the ``${4\pi}$-theorem'',
due to J.C.C. Nitsche \cite{Ni}.

\begin{theorem}[Nitsche] Let $\gamma$ be a real analytic curve in 
  ${\R}^3$ having total curvature less than $4\pi$.  Then $\gamma$
  bounds a unique minimal surface which is topologically equivalent to
  a disk, and this disk is stable and has no branch points.
\end{theorem}

\vskip .6cm

{\em Sketch of Proof.}  It follows easily from previous work of Shiffman
  that if $\gamma$ bounded more than one minimal disk which is a
  $C^0$-minimum to area, then it would also bound a branched minimal
  disk that is not a $C^0$-minimum to area.  Nitsche proved, for the
  curve $\gamma$ stated in the theorem, that every minimal disk which
  it bounds is a $C^0$-minimum for area, and hence by Shiffman's
  result it must be unique.\ \ \ \ \ \ \ \ \ \ \ \ \ \ \ \ \ \ \ \ \ \ 
  \ \ \ \ \ \ \ \ \ \ \ \ \ \ \ \ \ \ \ \ \ \ \ \ \ \ \ \ \ \ \ \ \ \ \ 
  \ \ \ \ \ \ \ \ \ \ \ \ \ \ \ \ \ \ \ \ \ \ \ \ \ \ \ \ \ \ \ \ \  q.e.d.

\setcounter{chapter}{3}
\setcounter{section}{-1}
\setcounter{theorem}{0}
\setcounter{lemma}{0}
\setcounter{definition}{0}
\setcounter{claim}{0}
\setcounter{corollary}{0}
\setcounter{condition}{0}
\setcounter{question}{0}
\setcounter{example}{0}
\setcounter{remark}{0}

\section{Existence and description}

In this section we show the existence of the minimal surfaces
${\Sigma}(V, t)$, and describe them geometrically, for sufficiently
small values of $t$.

\section{Existence}

Recall that $P_t$ denotes the plane ${\{} (x, y, z){\in}{{\R}^3} | z =
t{\}}$.  Let us consider two closed smooth Jordan curves, ${\alpha}$
and ${\beta}={\beta}(0)$ in the plane ${P_0}:z=0$ of ${\R}^3$.
Suppose that the curves ${\alpha}$ and ${\beta}$ intersect
transversely in a finite number of points $p_1,\ p_2, {\ldots},\ p_N$,
which we will call {\it crossing points}. 

For our purposes in this paper, the word {\em varifold} will denote a
finite union of connected compact planar regions which are smooth
except possibly at a finite number of boundary points, and to the
interior of each such component we associate a nonnegative integer,
called the {\em multiplicity} of that component, without orientation.
We will denote by ${\beta}(t)$ the translated image of $\beta$ by the
vector $(0, 0, t)$.

First of all we will prove an existence theorem, which shows that the
class of minimal surfaces that are our object of study is not the
empty one, and that it makes sense to study the asymptotic behavior.

\begin{theorem} Let $\alpha$ and
 ${\beta}$ be two Jordan curves in the
  plane $P_0$.  Let $V$ be a fixed varifold in $P_0$ with
  ${\Z}_2$-boundary given by $Z = {\alpha}\cup{\beta}$.  Then, for
  every $t$ sufficiently close to $0$, there exists a compact stable
  embedded minimal surface ${\Sigma}(t)$ whose boundary is
  ${\alpha}{\cup}{\beta}(t)$, and such that $\Sigma(t)$ is ``close
  to'' $V$, in the sense that the sequence of minimal surfaces
  ${\{}{\Sigma}(t){\}}_t$ converges to $V$, as $t {\to} 0$.
\end{theorem}

The proof will follow from the following theorem proven by Meeks and
Yau \cite{MeYa1}, which we state here without proof.

\begin{theorem}[Meeks-Yau] Let $M$ be a submanifold of a
$3$-dimensional analytic manifold whose boundary is piecewise smooth
and has nonnegative mean curvature with respect to the inner pointing
normal vector, and has interior angles less than $\pi$.  Let $\Sigma$
be a compact subdomain of $\partial M$ such that $\Sigma$ is
incompressible in $M$, namely each homotopically nontrivial curve in
$\Sigma$ is also homotopically nontrivial in $M$.  Then there exists
a stable minimal embedding $f:{\Sigma}\to M$ such that
$f(\partial{\Sigma}) = \partial{\Sigma}$.  Moreover if
$g:{{\Sigma}'}\to M$ is another minimal immersion of a compact surface
such that $g(\partial{{\Sigma}'})=\partial{\Sigma}$, then one can
assume $f(\Sigma){\subset}M{\backslash}g({\Sigma}')$.
\end{theorem}

{\em Proof of Theorem 3.1.}  Let ${\cal R}(\alpha)$ and
  ${\cal R}(\beta)$ be the bounded simply connected components of
  ${P_0}{\backslash}{\alpha}$ and ${P_0}{\backslash}{\beta}$
  respectively, and let $\cal R$ be the region in the plane $P_0$
  given by $\cal R = {\cal R}(\alpha){\cup}{\cal R}(\beta)$. 
  
  We wish to construct the analytic three-manifold $M(V, t)$
  satisfying the hypotheses of theorem $3.2$, for $t$ sufficiently
  small.  Let $C_1, \ldots , C_k$ be the bounded components of
  ${P_0}{\backslash}Z$ with $V$-multiplicity equal to zero.  Choose
  round circles ${S_i}\subseteq{Int}(C_i)$, for each $i{\in}{\{}1,
  \ldots , k{\}}$, and let $S_i(t)$ be the vertical upward translated
  of $S_i$ by $t$.  Since the number of connected components of $P_0
  {\backslash} ({\alpha}{\cup}{\beta})$ is finite, there will be a
  largest $t'$ such that for $t {\leq} t'$ there exist stable
  catenoids $A_i(t)$ bounded by ${S_i}{\cup}{S_i}(t)$ whose height is
  $t$.  We do not place catenoids above any component of $V$ which has 
  multiplicity $2$, nor above components of $V$ whose multiplicity 
  is $1$.
  
  Next, glue flat components away from the $(0, 1, 0, 1)$ crossing
  points, and where the prescribed multiplicities are $1$ or $2$.
  The glueing above the multiplicity $2$ components should be made
  so that the following condition is satisfied.
  
  Let $\sigma$ be the multiplicity two arc bounded by the two 
  $(0, 1, 2, 0)$ crossing points.  Then of the two flat pieces 
  glued above the multiplicity $2$ component, one does not contain 
  $\sigma$, and it glues on smoothly (as its projection ``crosses'' 
  $\sigma$) to the adjacent flat piece which also doesn't contain 
  $\sigma$ and projects on the multiplicity $1$ component ``across'' 
  $\sigma$.
  
  Finally, along the $(0, 1, 0, 1)$ crossing points, glue the corresponding
  multiplicity $1$ pieces via a twisted strip (the glueing taking place along
  the ``helices''), in such a way that if $\cal S$ is the surface just
  constructed, then the curvature of the strip with respect to the
  inner pointing normal vector (with respect to $\cal S$) is nonnegative.
  
  Then the three-manifold $M = M(V, t)$ is given by

  $$M = 
      Slab(t){\backslash}
            ({\bigcup}_{i=1}^{k}Int(A_i(t)){\bigcup}(Out({\cal S}))),
  $$
  where $Slab(t)$ is the region of space between the planes $P_0$ and
  $P_t$, $Int(A_i(t))$ is the bounded simply connected region in
  $Slab(t)$ containing $A_i(t)$ in its boundary, $Out(\cal S)$ is the
  unbounded connected component of $Slab(t){\backslash}{\cal S}$.
  Note that $M$ is analytic because the boundaries of the catenoids
  are circles.  Also, by construction, the interior angles are less
  than $\pi$, and the curvature with respect to the inner pointing
  normal is nonnegative.

  Then the surface $\Sigma$ satisfying the hypotheses of
  theorem $3.2$ is going to be $\cal S$.

  By construction, the surface $\Sigma$ constructed in this way has

  $${\pi}_{1}(\Sigma) = {\cal F}(k),$$ where ${\cal F}(k)$
  denotes the free group on $k$ generators.  Moreover $\Sigma$ admits a
  submersion onto a planar region in the plane $P_0$ which is
  homotopically equivalent to $V{\subseteq} \cal R$.  Now, $V$ is
  easily seen to be homotopically equivalent to $M$ (which is,
  homotopically, the cross product of the unit interval with a plane
  from which $k$ disks have been removed), and so one has

  $${\pi}_{1}(M) = {\cal F}(k).$$This allows us to conclude that $\Sigma$ is
  incompressible in $M$ because of the commutativity of the following
  diagram:

\begin{picture}(150,100)(0,0)
\put(45, 80){${\pi_1}(\Sigma)$}
\put(150,80){${\pi_1}(M)$}
\put(97,10){${\pi_1}(V)$}

\put(55,45){1-1}
\put(59,72){\vector(1,-1){48}}
\put(161,72){\vector(-1,-1){48}}
\put(161,72){\vector(1,1){1}}
\put(152,45){$isomorphism$}
\put(76,84){\vector(1,0){70}}
\put(110,60){\oval(10,10)[tr]}
\put(110,60){\oval(10,10)[tl]}
\put(110,60){\oval(10,10)[br]}
\put(110,55){\vector(-1,0){1}}
\end{picture}

  Hence the above construction satisfies all the hypotheses of Meeks
  and Yau's theorem, providing us with a {\it good barrier for the
    solution of Plateau's problem}.  We are assured the existence of a
  compact stable embedded minimal surfaces ${\Sigma}(t)$, having
  boundary given by ${\alpha}{\cup}{\beta}(t)$, for all $t {\leq}t'$,
  and incompressible in $M$.  By construction, we have that ${lim}_{t
    \to 0}({\Sigma}(t)) = V$.

{\hskip 11cm} q.e.d.

\begin{remark} As the number of annuli used as barriers in the previous 
  proof increases, the genus of the minimal surfaces ${\Sigma}(V, t)$
  correspondingly obtained increases, and the area of ${\Sigma}(V, t)$
  decreases, since the area of $V$ decreases.  Hence the least area
  surfaces bounded by ${\alpha}{\cup}{\beta}(t)$ are the ones having
  largest genus, and they correspond to the only varifold which has no
  components with multiplicity $2$.  The surfaces of largest area
  bounded by ${\alpha}{\cup}{\beta}(t)$ correspond to the disjoint
  union of the regions ${\cal R}(\alpha)$ and ${\cal R}({\beta}(t))$.
\end{remark}

\section{Description}

In this section we will describe geometrically the minimal
surfaces whose existence has been shown in the previous section, when
the distance between the two planes containing the boundary of such
surfaces is sufficiently small.

Let's give here a definition which will be useful in the proof of the
next theorem.

\begin{definition} A sequence of surfaces $\{{S_i}\}_{i=1}^{\infty}$
  with boundary is said to converge to a proper surface $S$ in
  ${\R}^3$ if, for each compact region $B$ of ${\R}^3$, there exists a
  positive integer $N_B$ such that, for $i>{N_B}$, ${S_i}{\cap}B$ is a
  normal graph over $S$, and $\{{S_i}\cap{B}\}_{i=1}^{\infty}$
  converges to $S{\cap}B$ in the $C^1$ norm
  (${\|}f{\|}_{C^1}={\|}f{\|}_0+{\|}Df{\|}_0$).  Moreover, if
  $\{S_i\}_{i=1}^{\infty}$ is a sequence of surfaces with boundary,
  and $\{{\phi_i}\}_{i=1}^{\infty}$ is a sequence of homotheties of
  ${\R}^3$ with center the origin, and
  $\{{\phi_i}(S_i)\}_{i=1}^{\infty}$ converges to a connected piece of
  a helicoid bounded by two straight lines, we will say that $S_i$ is
  {\it approximately helicoidal}\ \ near the origin for $i$
  sufficiently large.  Similarly one can define the convergence of a
  continuous family $\{S_t\}$, as $t$ approaches $t_0$.
\end{definition}

The minimal surfaces considered in the Main Theorem were
partially described by Wayne Rossman \cite{Ro}, under the hypotheses
that these minimal surfaces exist and are actually area-minimizing.
By theorem $3.1$ we have proven the existence of stable embedded
minimal surfaces which are not necessarily area-minimizing (see the above 
remark).  In this section we prove theorem 3.7, which will provide a
complete description of the {\it stable} minimal surfaces under
consideration in this paper.

Before stating the theorem, let us make a few observations about the
crossing points.  More precisely, consider a varifold $V$, arising as
the limit of a family ${\Sigma}(t_i)$ of minimal surfaces, in a small
neighborhood of the crossing points.  This is pictured schematically
in Figure $3.1$, where $a, b, c, d$ represent the multiplicity of $V$ in
each of the four regions.  Then the following claim holds.

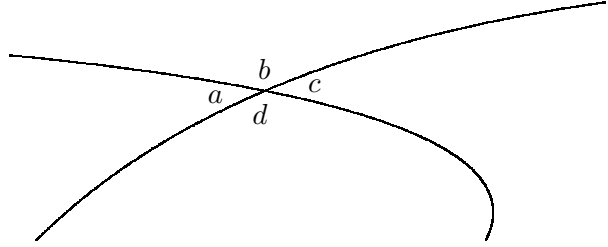
\begin{figure}
\begin{picture}(250,100)(0,0)
\qbezier(20,80)(225,60)(200,10)
\qbezier(30,10)(100,80)(245,100)
\put(95,62){$a$}
\put(114,71){$b$}
\put(133,66){$c$}
\put(112,54){$d$}
\end{picture}
    \caption{The varifold $V$, with multiplicity, around a crossing point.}
 \end{figure}

\begin{claim}  Let $V$ be a varifold determined by any of the minimal 
surfaces obtained in theorem $3.1$.  Then, with notation as introduced
above, the following three properties hold for $V$: 

\begin{description}
\item{-} The difference between the multiplicity of any region and an
  adjacent one in Figure $3.1$ is exactly one.

\item{-} At least one of $a, b, c, d$ must be equal to zero.

\item{-} The only possibilities allowed for the multiplicity at the
crossing points are $(0, 1, 0, 1)$ and $(0, 1, 2, 1)$.
\end{description}

\end{claim}

\begin{proof} The proof of the three assertions in the claim follows
  easily from the fact that the surfaces being considered are embedded
  and stable, and that $Z$ is the homology boundary of these surfaces.
\end{proof}

Hence we will, from now on, speak of points with multiplicity $(0, 1,
0, 1)$ and $(0, 1, 2, 1)$.  The two cases are illustrated in Figure
$3.2$.  

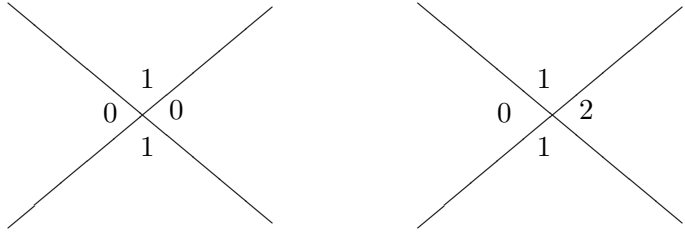
\begin{figure}
\begin{picture}(250,100)
\put(30,90){\line(6,-5){100}}
\put(185,90){\line(6,-5){100}}
\put(30,5){\line(6,5){100}}
\put(185,5){\line(6,5){100}}
\put(66,45){0}
\put(80,58){1}
\put(91,46){0}
\put(80,32){1}

\put(215,45){0}
\put(230,58){1}
\put(246,46){2}
\put(230,32){1}

\end{picture}

   \caption{Crossing points with multiplicity $(0, 1, 0, 1)$ and 
     $(0, 1, 2, 1)$.}
\end{figure}

\begin{remark} The previous claim implies that, for $t$ sufficiently 
  small, the varifold associated with each embedded compact stable minimal
  surface \newline bounded by ${\alpha}{\cup}{\beta}(t)$ is actually in the set
  ${\cal V}({\alpha}, {\beta})$, whose meaning is the same as in
  the introduction.\end{remark}

\begin{theorem}  Let ${\Sigma}(t)$ be a compact stable embedded minimal
  surface with boundary given by the union of the curves ${\alpha}$
  and ${\beta}(t)$.  Then, for $t$ sufficiently close to $0$, the
  surface ${\Sigma}(t)$ has the following properties:

\begin{description}
\item{1.} In the complement of the union of small cylindrical
  neighborhoods $N(i)$ of the crossing points, the connected
  components of ${\Sigma}(t){\backslash}{\cup}N(i)$ are graphs over
  their projection to the plane $P_0$.

\item{2.} In $N(i)$, ${\Sigma}(t)$ is either approximately helicoidal,
  or it is the union of two graphs over the plane $P_0$.
\end{description}
\end{theorem}

\begin{proof}  First notice that the theorem certainly holds if 
  ${\Sigma}(t)$ consists of the union of the two planar regions
  bounded by ${\alpha}$ and ${\beta}(t)$.  Hence we need to prove the
  theorem for connected ${\Sigma}(t)$.  The proof will be divided in
  five steps.

  \underline{Step 1}: {\it For $t$ sufficiently close to $0$, the
    components of ${\Sigma}(t)$ away from the crossing points are
    graphs over $P_0$.}

For each $1{\leq}i{\leq}n$, let $D_{i}(\epsilon)$ be a disk in $P_0$
with center $p_i$ and radius $\epsilon$, chosen so that the following
two conditions are satisfied:

\begin{description}
\item(a) the disks $D_i$ are mutually disjoint;

\item(b) for each $i$, ${D_i}{\cap}({\alpha}{\cup}{\beta})$ is
  topologically a ${\times}$.

\end{description} Let 
$N={P_0}{\backslash}{\bigcup}_{i=1}^{N}{D_i(\epsilon)}$, and let 
$U=N{\times}{\R}$ be the vertical cylinder over $N$.  Moreover, let 
${\hat {\Sigma}}(t)={{\Sigma}(t)}{\cap}U$.  Let us notice that along 
the curves ${\alpha}{\cap}{\hat {\Sigma}(t)}$ and 
${\beta}(t){\cap}{\hat {\Sigma}(t)}$, the normal vector to 
${\Sigma}(t)$ must become arbitrarily vertical as $t$ approaches $0$.  
To see this, for example along ${\beta}(t){\cap}{\hat {\Sigma}(t)}$, 
suppose $\Sigma(t)$ is vertical at a boundary point $p$ of 
${\beta}(t){\cap}{\hat {\Sigma}(t)}$.  Consider, above and away from 
${\Sigma}(t)$, a piece of a half strictly unstable vertical catenoid 
with boundary curves a larger circle in the plane $P_0$ and a smaller 
circle in the plane $P_{t+{\epsilon}'}$ (this catenoid is a graph 
above the plane $P_0$, except along the smaller boundary circle).  
Translate vertically down this catenoid piece until it makes first 
contact with ${\Sigma}(t)$ at some point $p$; the maximum principle 
implies first that $p$ is a boundary point both of the catenoid piece 
and of ${\Sigma}(t)$, and then that the two minimal surfaces 
${\Sigma}(t)$ and the catenoid piece coincide in a small neighborhood 
of $p$.  This contradiction shows that $\Sigma(t)$ cannot become 
vertical along its boundary, away from the crossing points.  Moreover 
we have that locally the surface $\hat {\Sigma}(t)$ is situated on the 
same side with respect to the vertical cylinder 
$(({\alpha}{\cup}{\beta}){\cap}{\hat {\Sigma}(t)}){\times}{\R}$, 
otherwise one could translate a catenoid in such a way that the first 
point of contact would be an interior point, which contradicts the 
maximum principle.  With the same notation as above, let us notice 
that there exists the upper lower bound of the angle $\theta(p,t)$ 
formed by ${\Sigma}(t)$ with an arbitrary catenoid piece intersecting 
${\Sigma}(t)$ only in $p$; in fact this infimum is given by the angle 
between the tangent plane to ${\Sigma}(t)$ at $p$ and the horizontal 
plane $P_0$.  Moreover, for an arbitrary $p$ in 
$({\alpha}{\cup}{\beta}(t)){\cap}{\hat {\Sigma}(t)}$, the upper lower 
bound of $\theta(p,t)$ approaches zero as $t$ gets closer to $0$.  Let 
$$\theta_0(t)={\max}{\{}\theta(p,t) |{\,} 
p{\in}({\alpha}{\cup}{\beta}(t)){\cap}{\hat {\Sigma}(t)}{\}}.$$ The 
compactness of $({\alpha}{\cup}{\beta}(t)){\cap}{\hat {\Sigma}(t)}$ 
guarantees that $\theta_0(t)$ is well defined.  Furthermore, what has 
been said above implies that $${\lim}_{t\to 0}{\theta_0}(t)=0,\ \ \ \ 
{\mbox{and therefore}}\ \ \ \ {\lim}_{t\to 0}{\tan}{\theta_0}(t)=0.$$ 
By R.  Schoen's estimate \cite{Sc2} there exists a universal constant 
$c$ such that $K{\!}<{\!}\frac{c}{r^2}$, where $K$ is the Gaussian 
curvature at a point of a stable minimal surface and $r$ is the 
distance between that point and the boundary of the surface.  Now let 
us choose $t$ sufficiently close to $0$ so that 
${\tan}{\theta_0}(t){\!}<{\!}{\min}{\{}{\frac{\pi}{32{\sqrt{c}}}},{\frac{1}{16}}{\}}$. 
 Let $p$ be a point in $\hat {\Sigma}(t)$, and let $\hat p$ be the 
orthogonal projection of $p$ on the plane $P_0$.  Let $r$ be the 
horizontal distance between $\hat p$ and 
${\alpha}{\cup}{\beta}{\backslash}{{\bigcup}_{i=1}^{n}{D_i}}$.  Let 
$Cyl_1$ be the part of the vertical cylinder $D({\hat 
p},{\frac{r}{2}}){\times}{\R}$ with height $2{r}{}{\tan}{\theta_0}$ 
and containing ${\hat {\Sigma}(t)}{\cap}{\partial}Cyl_1$ in its 
interior.  Let $Cyl_2$ be the part of the vertical cylinder $D({\hat 
p},{\frac{r}{4}}){\times}{\R}$ with height $2{r}{}{\tan}{\theta_0}$ 
which is contained in $Cyl_1$.  Let us suppose that our assertion in 
Step 1 is not true.  It will be therefore possible to find a point $q$ 
in ${\Sigma}(t){\cap}Cyl_2$ whose normal vector ${\vec N}(q)$ is 
horizontal, namely $<{\vec N(q)}, {\vec {e_3}}>=0$.  There follows the 
existence of a geodesic ${\omega}(s)$ in ${\Sigma}(t)$ with unimodular 
velocity, and such that ${\omega}(0)=q$, ${\omega}'(0)={\vec e_3}$, 
and the curvature $|{\omega}''(s)|$ of ${\omega}(s)$ in $Cyl_1$ is 
bounded above by $\frac{2{\sqrt{c}}}{r}$ because of Schoen's estimate.


Notice that, since ${\tan}{\theta_0}<\frac{1}{16}$ (and since
$|{\omega}'(s)|=1$ implies $|d{\omega}|{\approx}|ds|$), ${\omega}(s)$
lies in $D({\hat p},{\frac{r}{2}}){\times}{\R}$ for
$0{\leq}s{\leq}{4}{r}{}{\tan}{\theta_0}$.  Considering the estimate
for the curvature of ${\omega}(s)$, the estimate
$\frac{\pi}{32{\sqrt{c}}}$ for ${\tan}(\theta_0)$, and the integral
$$\int_{0}^{4{r}{}{\tan}{\theta_0}}|{\omega}''(s)|{\,}ds,$$ we can
conclude that the length of the curve ${\omega}'(s)$ in the unit
sphere $S^2$ is less than $\frac{\pi}{4}$, for
$0{\leq}s{\leq}4{r}{}{\tan}{\theta_0}$.  Since ${\omega}'(0)={\vec
  e_3}$, we have that $<{\,}{\omega}'(s),{\vec e_1}{\,}>$ and
$<{\,}{\omega}'(s),{\vec e_2}{\,}>$ are both less than $\frac{\sqrt
  2}{2}$, and that $<{\,}{\omega}'(s),{\vec e_3}{\,}>$ is larger than
$\frac{\sqrt 2}{2}$, for $0{\leq}s{\leq}4{r}{}{\tan}{\theta_0}$.  Hence it is
$${\langle}{\omega(4{r}{\tan}{\theta_0})}-{\omega}(0),{\vec{e_3}}{\rangle}=
\int_{0}^{4{r}{\tan}{\theta_0}}{\langle}{\omega}'(s),{\vec
  e_3}{\rangle}ds> {\frac{\sqrt
    2}{2}}4{r}{}{\tan}{\theta_0}>2{r}{\tan}{\theta_0}.$$ Likewise we
have:
$$|{\langle}{\omega}(4{r}{}{\tan}{\theta_0})-{\omega(0)},{\vec
  e_1}{\rangle}|< \frac{r{\sqrt 2}}{8},$$ and
$$|{\langle}{\omega}(4{r}{}{\tan}{\theta_0})-{\omega(0)},{\vec
  e_2}{\rangle}|< \frac{r{\sqrt 2}}{8},$$ because of the condition
${\tan}{\theta_0}{\!}<{\!}\frac{1}{16}$.  But this implies that
$\omega(s)$ intersects one of the two horizontal disks of
${\partial}Cyl_1$ for some value of $s$ between $0$ and
$4{r}{}{\tan}{\theta_0}$, which is impossible, given the way $Cyl_1$
was constructed.  Hence $\vec N$ is never horizontal on $\hat {\Sigma}(t)$,
which therefore is union of graphs over $D({\hat p},\frac{r}{4})$.

Moreover, notice that at most two sheets can lie over $D({\hat
  p},\frac{r}{4})$, otherwise ${\Sigma}(t)$ would not be embedded (at
least along the boundary).  If the surface is area-minimizing, then at
most one sheet lies above each point, but if no assumption of area
minimality is made, then there could be points above which there are
two sheets.

Let us remark that for $t\approx 0$ the above graphs are ``almost
horizontal'', since it is clear the above argument can be strengthened
to show that, for any $x {\in} (0, 1)$, it is possible to make
$|{\langle}\vec N{(q)}, \vec{e_3}{\rangle}| > x$ for each $q {\in}
{\Sigma}(t)$, for $t$ sufficiently close to $0$.

\underline{Step 2}: {\it ${\Sigma}(t)$ is the union of two disjoint
  graphs in the neighborhood of each $p$ in
  ${\partial}{V_1}{\cap}{\partial}{V_2}$, where $V_1$ and $V_2$ are
  components of $V$ with multiplicities $1$ and $2$ respectively.}

To see this, consider a small vertical cylinder $\cal C$ whose 
vertical axis contains $p$ and whose height is $t$.  Moreover, let 
$\hat {\Sigma}(t) = {\Sigma}(t){\cap}{\cal C}$.  Now homotetically 
expand ${\cal C}$, with coefficient of homothety $1/t$, and center of 
the homothety in $p$.  The image of $\hat {\Sigma}(t)$ under the 
homothety converges, as $t{\to}\infty$, to a minimal surface having a 
simply connected component $S$ which is a stable minimal surface 
contained in a half-space and bounded by a straight line, $\ell$.  The 
image $g(\ell)$ of $\ell$ via the Gauss map is either a single point 
or a great-circle on $S^2$ containing the north and south pole of 
$S^2$, and near the boundary $\ell$ the image of $S$ via the Gauss map 
is entirely contained in one of the two hemispheres determined by 
$g(\ell)$, because of the way the barrier to get ${\Sigma}(t)$ was 
constructed (theorem $3.1$).  Hence the image under the Gauss map of 
$S$ is entirely contained in such half hemisphere, by the hypothesis 
of stability.  By reflection with respect to the line $\ell$ one 
obtains a complete minimal surface containing a line and having total 
curvature between $- 4{\pi}$ and $0$.  The two only complete minimal 
surfaces with total curvature $-4{\pi}$ are the catenoid and Enneper 
surface \cite[p.40]{BaCo}.  Since the catenoid does not contain a line of 
reflective symmetry, and Enneper surface's Gauss map does not satisfy
our conditions, then $S$ must be a half plane, and its image via 
the Gauss map must ba a point.  Finally, notice that the points 
contained in the other simply connected component of the homothetic 
expansion of $\hat {\Sigma}(t)$ correspond to points which are 
contained in the interior of ${\Sigma}(t)$, and hence Schoen's 
curvature estimate applies to them.  These observations complete the 
proof.

\underline{Step 3}: {\it ${\Sigma}(t)$ is the union of two disjoint
  graphs in the neighborhood of each crossing point having
  multiplicity $(0, 1, 2, 1)$.}

This follows from the previous two steps, and from the analysis of the
possible liftings to ${\Sigma}(t)$ of small circles in $V$ around the
point under consideration.


Our next aim is to study ${\Sigma}(t)$ around a crossing point whose
multiplicity is $(0, 1, 0, 1)$.

\underline{Step 4}: {\it ${\Sigma}(t)$ is topologically a disk in the
  neighborhood of each crossing point with multiplicity $(0, 1, 0,
  1)$.}

We will show that in a closed cylindrical neighborhood $U$ of each
crossing point with multiplicity $(0, 1, 0, 1)$, the compact surface
${\Sigma}(t){\cap}U$ has genus zero for $t$ sufficiently close to $0$.
In order to do this, let us homothetically expand the spherical
neighborhood $U$ with center in a point $p$, with expansion
coefficient $\frac{1}{t}$, where $p$ is a point of maximum Gaussian
curvature inside ${\Sigma}(t){\cap}U$.  Let us denote by $\tilde U$
the expanded neighborhood, and notice that the expansion transforms
the planes $P_0$ and $P_t$ to two new planes which have distance equal
to $1$ from each other.  Moreover the homothety takes the arcs
${\alpha}{\cap}U$ and ${\beta}(t){\cap}U$ to arcs $\tilde {\alpha}$
and $\tilde{\beta}(t)$ in $\tilde U$ which are segments of an almost
straight line.  Let ${\partial}{\tilde{\alpha}} = {\{} a_1, a_2 {\}}$,
and ${\partial}{\tilde{\beta}(t)} = {\{} b_1, b_2 {\}}$, with the
convention that the orthogonal projections of $a_1$ and $b_1$ on $P_0$
lie in the boundary of the same multiplicity one component of the
interior of $V$, and the orthogonal projections of $a_2$ and $b_2$ lie
in a different multiplicity one component.  Now let us join $a_1$ to
$b_1$ and $a_2$ to $b_2$ by two geodesic arcs contained in the
homothetic expansion $\tilde {\Sigma}(t)$ of ${\Sigma}(t)$.  Because
${\Sigma}(t)$ is a graph over $P_0$ away from the crossing points with
multiplicity $(0, 1, 0, 1)$, and because for each $x{\in}(0,1)$ we
have $|{\langle}{\vec N}(q), {\vec e}_3{\rangle}|>x$ away from these
crossing points, for $t{\approx}0$ (namely the normal vector to
${\Sigma}(t)$ in $q$ is almost vertical if $q$ is away from the
crossing points with multiplicity $(0, 1, 0, 1)$), we may assume that
the two geodesic arcs defined above project orthogonally on two
different multiplicity one components of the interior of $V$, and
hence do not intersect each other.  We will now show that the piece
$\bar {\Sigma}(t)$ contained in $\tilde {\Sigma}(t)$ and bounded by
the loop union of the four curves $\tilde {\alpha}$, $\tilde
{\beta}(t)$, and the two geodesic arcs previously defined, is a disk.
Since $\partial({\Sigma}(t))$ is contained inside the boundary of the
convex hull of ${\Sigma}(t)$, ${\Sigma}(t)$ separates its convex hull
(which is simply connected) into two distinct types of regions, one
associated with the ``$+$'' sign and the other with the ``$-$'' sign.
Hence ${\Sigma}(t)$ is orientable, and consequently also 
$\bar{\Sigma}(t)$ is so.  This implies that the Gauss map, from the 
oriented $\bar{\Sigma}(t)$ to the unit sphere $S^2$, is well defined.
The Gaussian image ${\nu}(t)$ of the boundary curve 
${\partial}\bar{\Sigma}(t)$ is a curve that lies in
a small neighborhood of the spherical region bounded by the union of
two great semicircles joining the north and south poles of $S^2$.  
Since ${\Sigma}(t)$ is stable,
$\bar{\Sigma}(t)$ is also stable.  Because the Gaussian image of
a stable minimal surface cannot contain a hemisphere, the
image of $\bar{\Sigma}(t)$ under the Gauss map can only contain one of the
two regions in the complement of this neighborhood in $S^2$, and must
be disjoint from the other region.  Since the winding number of the
Gauss map around ${\partial}\bar{\Sigma}(t)$ is one, the Gauss map can only
cover this region once.  This implies that for $t$ sufficiently close
to $0$ it is:
$$-2{\pi} {\leq} {\int}_{\bar{\Sigma}(t)}{K}{\,}dA < 0.$$ Now, the geodesic
curvature $k_g$ is zero along the two almost straight arcs contained
in $\partial \bar{\Sigma}(t)$, and along $\tilde {\alpha}$ and
$\tilde{\beta}(t)$ the geodesic curvature is approximately equal to
zero.  The sum of the exterior angles where these smooth arcs
intersect is between $0$ and $4{\pi}$.  Therefore by the Gauss-Bonnet
theorem, one has that the Euler characteristic of $\bar{\Sigma}(t)$ is
either zero or one.  Since ${\partial}\bar{\Sigma}(t)$ consists exactly of
one curve, we can conclude that the Euler characteristic is $1$, and
that $\bar{\Sigma}(t)$ is topologically a disk.  Hence, in small
neighborhoods of the crossing points with multiplicity $(0, 1, 0, 1)$,
${\Sigma}(t)$ is topologically a disk.

\underline{Step 5}: {\it ${\Sigma}(t)$ is approximately helicoidal around
  crossing points having multiplicity $(0, 1, 0, 1)$.}

Normalize ${\Sigma}(t){\cap}N$ by a homothety with center a
point $p$ of maximum Gaussian curvature, in such a way that
${\max}_{q{\in}{\Sigma}(t){\cap}N}\{|K(q)|\}=1$ on the normalized
surface, which we shall denote by ${\breve {\Sigma}}(t)$.  Modulo a
translation, we can suppose that the point $p{\in}\beta$ is the
origin.  Let us notice that
${\max}_{q\in{{\Sigma}(t)}{\cap}N}\{|K(q)|\}\to{\infty}$ on
${\Sigma}(t){\cap}N$ as $t$ approaches $0$, since the Gauss map $\nu$
is almost vertical on $\partial{{\Sigma}(t)}$, except near the
crossing points where $\nu$ changes very quickly, and hence the
modulus of $K={Jac}(d{\nu})$ must be large near the crossing points.
Therefore normalizing ${\max}_{q\in{{\Sigma}(t)}{\cap}N}\{|K(q)|\}=1$
involves a dilation factor which becomes arbitrarily large as $t$
approaches $0$, and hence the two planar curves in the boundary of
${\breve {\Sigma}}(t)$ become arbitrarily straight as $t$ approaches
$0$, around the crossing points.  A result of Anderson \cite{An}
states that for each sequence of surfaces
${\{}{\breve{{\Sigma}}}_{t_{\ell}}{\}}_{{\ell}=1}^{\infty}${},
$t\to{\infty}$, it is possible to extract a convergent subsequence (in
the $C^2$ norm) $\{S_{ij}\}_{i=1}^{\infty}$ in the compact spherical
neighborhood $B(0,j)$ with radius $j$ in ${\R}^3$.  This sequence of
surfaces can be chosen in such a way that ${\{}S_{ij}{\}}_{i =
  1}^{\infty}$ is a subsequence of ${\{}S_{kl}{\}}_{k = 1}^{\infty}$
if $j{\,}>{\,}l$.  The sequence ${\{}S_{mm}{\}}_{m = 1}^{\infty}$ is
the Cantor diagonalization, and converges in the $C^2$ norm in
arbitrary compact regions to a surface $S$ having one or two boundary
curves, which must be straight lines.  Moreover Anderson's compactness
theorem \cite{An} implies that $S$ is embedded.  The limit surface $S$
is simply connected in any compact spherical neighborhood, and
therefore is simply connected in ${\R}^3$.  If the boundary of $S$
consists of two straight lines, then by Alexandrov reflection
principle $S$ can be extended to a simply connected minimal surface in
${\R}^3$ properly embedded, without boundary and with infinite
symmetry group.  By virtue of a theorem of Meeks and Rosenberg
\cite{MeRo}, such extended minimal surface is a plane or a helicoid.
However this surface contains two straight lines which do not
intersect and are not parallel to each other, and hence $S$ is a piece
of a helicoid, with one or two boundary lines.  If $S$ had only one
boundary line, then it could be extended via a rotation of angle $\pi$
around the straight boundary line, producing a properly embedded
minimal surface.  Because the convergence of the above subsequence to
$S$ is with respect to the $C^2$ norm, the normal vectors are
converging as well, and the extended surface has finite total
curvature.  By a result of L{\'{o}}pez and Ros such a surface must be
a plane or a catenoid; however since it is simply connected, it must
be a plane.  Hence $S$ is a half-space.  Since the convergence to $S$
is of class $C^2$, we know that the normal vectors along $\partial S$
are not constant as $m\to {\infty}$, which implies that $S$ cannot be
a plane.  Therefore $S$ must be a piece of helicoid with two boundary
lines.  Since ${\hat {\Sigma}(t)}$ is a graph over $P_0$ and because
we are considering a point with multiplicity $(0, 1, 0, 1)$, ${\hat
  {\Sigma}(t)}$ is a graph over the components of
${P_0}{\backslash}({\beta}(0){\cup}{\alpha})$ having the ``$+$'' sign.
Hence $S$ is totally determined.  If there was some subsequence which
would not be eventually contained in a given ${\epsilon}$-neighborhood
of $S$ in a given sphere $B(0,j)$ in ${\R}^3$, then it would be
possible to find a subsequence ${\{}{\check
  {\Sigma}_{t_{\ell}}}{\}}_{{\ell}=1}^{\infty}$ converging to a point
not belonging to $S$, which is a contradiction.  This allows us to
conclude that any sequence ${\{}{\check {\Sigma}_{t_{\ell}}}{\}}_{\ell
  =1}^{\infty}$ such that $t_{\ell}\to 1$ eventually lies in a
predetermined arbitrarily small $\epsilon$-neighborhood of $S$ inside
each compact region of ${\R}^3$.  For an arbitrarily given $\epsilon$,
let us choose $\ell$ big enough in such a way that the surface
${\check {\Sigma}_{t_{\ell}}}{\cap}{D_{j}(0)}$ be contained in an
$\epsilon$-neighborhood of $S{\cap}{D_{j}(0)}$.  For each pair of
points $p$ and $q$, with $p$ in $\check {\Sigma}_{t_{\ell}}$ with
normal vector $\vec N(p)$ to $\check {\Sigma}_{t_{\ell}}$ and $q$ in
$S$ with normal vector $\vec N(q)$ to $S$ such that
${dist}(p,q)<\epsilon$, the estimate on the function $|K(q)|$ implies
that $|<{\vec N(p)},{\vec N(q)}>|$ is bounded away from zero.  In fact
by choosing $\epsilon$ sufficiently close to zero and $\ell$
sufficiently large, we will be able to achieve $|<{\vec N(p)},{\vec
  N(q)}>|$ to be arbitrarily close to $1$.  It follows that $\check
{\Sigma}_{t_{\ell}}$ is union of graphs on $S$ for $\ell$ sufficiently
large, and hence $\check {\Sigma}_{t_{\ell}}$ is a one-sheeted graph
on $S$ around a crossing point with multiplicity $(1,0,1,0)$, which is
in accordance to the end of the proof of Part $1$.  This allows us to
conclude that ${\{}{\check {\Sigma}_{t_{\ell}}}{\}}_{\ell
  =1}^{\infty}$ converges to $S$ in the $C^0$-norm as one-sheeted
graphs, and that the normal vectors are convergent as well; hence we
have in addition that ${\{}{\check {\Sigma}_{t_{\ell}}}{\}}_{\ell
  =1}^{\infty}$ converges to $S$ in the $C^1$-norm in any compact
sphere, and that ${\Sigma}(t)$ is approximately helicoidal in a
neighborhood of each crossing point with multiplicity $(1,0,1,0)$, for
$t$ sufficiently close to $0$.

\end{proof}

\begin{question}  Does there exist, for $t$ sufficiently small, an
  unstable embedded minimal surface bounded by
  ${\alpha}{\cup}{\beta}(t)$, and having genus larger than any compact
  stable embedded minimal surface bounded by
  ${\alpha}{\cup}{\beta}(t)$?
\end{question}

The previous theorem described the surfaces ${\Sigma}(t)$
geometrically.  The next theorem describes their topological
properties.  Let $v_1$, $v_2$, $e_1$, $e_2$, $f_1$ and $f_2$ 
be defined as in section $1$.

\begin{theorem} Let $V{\in}{\cal V}({\cal A}, {\cal B})$, and for $t$
sufficiently small let ${\Sigma}(V, t)$ be a minimal surface given in
the statement of theorem $3.3$.  Then the Euler characteristic of
${\Sigma}(V, t)$ is equal to $({v_1} + 2{v_2}) - ({e_1} + 2{e_2}) +
({f_1} + 2{f_2}) $.  Since there is exactly one varifold
${V_0}{\in}{\cal V}({\cal A}, {\cal B})$ with $v_2 = e_2 = f_2 = 0$,
there is exactly one topological type ${\Sigma}({V_0}, t)$ in ${\cal
S}(t)$ with Euler characteristic equal to $v_1 - e_1 + f_1$.
\end{theorem}

\begin{proof} This follows easily from the description of ${\Sigma}(t)$, 
  and from the observation that a cell decomposition of ${\Sigma}(t)$
  can be computed via a cell decomposition of the unique varifold $V$
  determined by ${\alpha}{\cup}{\beta}$, and corresponding to
  ${\Sigma}(t)$.  Such a cell decomposition is equivalent to one that
  has:

\begin{description}
\item{---} number of vertices equal to $v_1 + 2v_2$.

\item{---} number of edges equal to $e_1 + 2e_2$.

\item{---} number of faces equal to $f_1 + f_2$.
\end{description}

The unique area-minimizing varifold gives rise to a minimal surface
having largest genus, since for this varifold $v_2 = e_2 = f_2 = 0$.
Hence the proof is finished.

\end{proof}

In sections $4$ and $5$ we will prove that the surface
${\Sigma}({V_0}, t)$ is the unique surface of least area in ${\cal
  S}(t)$, for $t$ sufficiently small.

\begin{remark} The results in this paper generalize
  to the case of boundary curves being two collections $\cal A$ and
  ${\cal B}(t)$ of disjoint smooth closed Jordan curves contained in
  the planes $P_0$ and $P_t$ respectively.  Moreover these results 
  also hold if $\cal A$ and ${\cal B}$ are contained in the interior
  of the half plane 
  ${P'}_0 = {\{} (x, y, z) \in {\R}^3 |\ \ z = 0, x{\geq}0{\}}$.
  In this case ${\cal B}(t)$ is the collection of curves in the plane
  ${P'}_t$, obtained by rotating $P_0$ counterclockwise around the
  $x$-axis by $t\frac{\pi}{2}$ radians.
\end{remark}

The above theorem provides a complete description of the stable
surfaces ${\Sigma}(t)$, description which we will use to prove the
uniqueness of the correspondence between ${\cal V}({\cal A}, {\cal B})$ 
and ${\cal S}(t)$ in the next two sections of this paper.

\setcounter{chapter}{4}
\setcounter{section}{-1}
\setcounter{theorem}{0}
\setcounter{lemma}{0}
\setcounter{definition}{0}
\setcounter{claim}{0}
\setcounter{corollary}{0}
\setcounter{condition}{0}
\setcounter{question}{0}
\setcounter{example}{0}
\setcounter{remark}{0}
\section{Proof of stability}

The main objective of this section is to prove the following Theorem.

\begin{theorem}  Let ${\Sigma}(t)$ be an embedded minimal surface whose 
  boundary is ${\alpha}{\cup}{\beta}(t)$.  Suppose that ${\Sigma}(t)$ is
  described as in Theorem $3.6$.  Then, for $t$ sufficiently close to
  $0$, ${\Sigma}(t)$ is stable.
\end{theorem}

\section{Preliminary lemmas}

The proof of theorem $4.1$ will follow from some preliminary lemmas.
To state these lemmas we will need some additional notation.  Let $U$
be a vertical cylindrical neighborhood with height larger than $1$ of
a helicoidal crossing point.  Let $E(t)$ be the intersection of $U$
with ${\Sigma}(t)$.  Suppose that the radius of $U$ is sufficiently
small, so that $E(t)$ is stable; note that this can always be
accomplished, for an appropriate choice of radius, because although
${\Sigma}(t)$ is approximately helicoidal in $U$, the boundary arcs of
the helicoidal piece are not parallel, and the multiplicity of the
crossing point contained in $U$ is $(0, 1, 0, 1)$; these conditions
guarantee that the area of the Gaussian image of $E(t)$, counted with
multiplicity, is less than $2\pi$.  Let $G(t)$ be one of the connected
components of $\hat{\Sigma}(t)$ (defined in Part $1$ of the proof of
Theorem $3.6$) adjacent to $E(t)$, and let ${\gamma}(t) = {\partial
E(t)}{\cap}{\partial G(t)}$.  Also recall that here $t$ approaches
$0$.  The lemma that we are about to prove guarantees that the
behavior of a supposed Jacobi vector field on ${\Sigma}(t)$ is not
``too wild'' away from the crossing points.  We remark here, once and
for all, that all the Jacobi vector fields we consider are
not identically zero, unless explicitly stated otherwise.

\begin{lemma}  Let ${\{}{\Sigma}(t){\}}_{t{\to}0}$ be a sequence of
  compact embedded minimal surfaces bounded by
  ${\alpha}{\cup}{\beta}(t)$, described as in Theorem $3.6$, and
  suppose that all $\Sigma(t)$ are unstable, with Jacobi vector fields
  given by ${u_t}:{\Sigma(t)}\to [0, {\infty})$, ${u_t}=0$ on
  ${\partial}{\Sigma}(t)$, $u_t$ of class $C^{\infty}$.  Then $u_t$ is
  not a ``bump function'' on ${\gamma}(t)$, for $t$ sufficiently close
  to $0$.  More precisely, if $p_t {\in} {\gamma}(t)$ is a local
  maximum for the function $u_t$, then $p_t$ is contained in an arc
  intersecting $E(t)$ along which the values attained by $u_t$ are
  very close to $u(p)$, of the order of  \newline  
  $1 - {\cos}({\vec N}(p_t), {\vec e}_3)$.
\end{lemma}

\begin{proof}  Since it is not ambiguous, after an appropriate choice
  of $t$, we will drop the $(t)$ indicating dependence on $t$, in this
  proof.  Let us consider the foliation ${\{}{\Sigma} +
  r{\vec{e_3}}{\}}_{0{\leq}r{\leq}1}$, and consider the surface
  ${\Sigma}'$ which is obtained by deforming $\Sigma$ via the Jacobi
  vector field ${\tilde u} = \frac{u}{{\max}u}$, namely the set of all
  the points $q+{\tilde u}(q){\vec N}(q)$, as $q$ varies in $\Sigma$.
  Let us consider the ``restrictions'' of the foliation and of the
  surface ${\Sigma}'$ to the subset $G$ of $\Sigma$, and let us denote
  by ${\{}G_r{\}}_{0 {\leq} r {\leq} 1}$ and $G'$ the transformed
  images of $G$ via the foliation and the Jacobi vector field,
  respectively.  Because all the surfaces we are considering are
  compact, there exists the maximum value $\overline r$ for which the
  surface $G'$ intersects the foliation ${\{}P_r{\}}_{0 {\leq} r
    {\leq}1}$ (after such value the translated surfaces
  ${\{}G_r{\}}_{r {\geq} {\overline r}}$ are situated ``above $G'$'').
  By the maximum principle at a boundary point, such last point of
  contact must be a point $\overline p$ contained in the boundary
  curve ${\gamma} = (\partial G){\cap}(\partial E)$.  Let $p$ be the
  point where ${\tilde u}_{|_{G}}$ attains its maximum value; since
  $G$ is an almost horizontal graph, there is no loss of generality in
  the assumption that the maximum value of ${\tilde u}_{|_{G}}$ is
  attained at $\overline p$ (in fact this can always be achieved by
  slightly deforming the curve $\gamma$).  Let now $S(\eta)$ be a very
  thin strip which is a ${\eta}$-neighborhood of $\gamma$ with respect
  to the metric of $\Sigma$, ${\eta}{\approx}0$, and $p'$ a point
  contained in $(\partial S(\eta){\cap}E){\backslash}{\partial
    {\Sigma}}$ corresponding to the last intersection point of
  ${\{}(P{\cup}S)_r{\}}_{0 {\leq} r {\leq} 1}$ with
  $(P{\cup}S(\eta))'$.  Let us notice that:
  $$
  1 - {\overline r} {\leq} 1 - {\cos}({\vec N}(p), {\vec e}_3),
  $$ and that certainly, if ${\overline r}'$ denotes the value of $0
  {\leq} r {\leq} 1$ corresponding to $p'$, then we have $1 -
  {\overline r}' {\leq} 1 - {\overline r}$, which implies that
  $$ |(p+{\tilde u}(p){\vec N}(p))_z - (p'+{\tilde u}(p'){\vec
    N}(p'))_z| {\leq} 1 - {\cos}({\vec N}(p), {\vec e}_3),
  $$ where $*_z$ denotes the ordinary $z$-coordinate of a point $*$ in
  ${\R}^3$.  Now, since the unit normal vector to $\Sigma$ approaches
  ${\pm}{\vec e}_3$ in a continuous fashion away from the crossing
  points, the above says actually that the difference between ${\tilde
    u}(p)$ and ${\tilde u}(p')$ is at most $1 - {\overline r}
  + |p_z - {p'}_z|$.  Finally, by perturbing $\gamma$ slightly in all
  directions around $p$ and applying the above argument to these
  perturbed curves, the proof of the lemma is complete.
\end{proof}

The next lemma is inspired by the paper \cite{Ka} of N. Kapouleas.  It
will show that if a Jacobi vector field on a manifold is ``close'' (in
some sense, which will be specified in the statement of the lemma) to
a Jacobi vector field on another appropriate manifold, then the first
eigenvalues of the Jacobi operator on the two manifolds are also
``close'' to each other.

\begin{lemma}  Let $U$ be a neighborhood of a helicoidal crossing point
  such that the surface $E(t)$ given by the intersection of ${\Sigma}(t)$
  with $U$ is stable, and let $G(t)$ be a connected component of
  $\hat{\Sigma}(t)$ adjacent to $E(t)$.  Let ${M_1} = {\Sigma}(t)$, and
  let ${M_2}(\mu)$ be ${{\Sigma}(t)}{\backslash}S(\mu)$, where $S(\mu)$
  is a $\mu$-neighborhood of $\gamma =
  {\partial}{E(t)}{\cap}{\partial}{G(t)}$ in the metric of ${\Sigma}(t)$,
  ${\mu}{\approx} 0$.  Let $f$ define a Jacobi vector field on $M_1$.
  Suppose that it is possible to deform the function
  $f{\colon}{M_1} \to {\R}$ to obtain a Jacobi vector field
  $G(f){\colon}{M_2}(\mu) \to {\R}$, which is zero on
  ${\partial}{M_2}(\mu){\cap}{\partial}{\Sigma}$ and such that it
  satisfies the three additional conditions:

\begin{description}

\item(i) ${\|}f{\|}_{\infty}{\leq}2{\|}G(f){\|}_{\infty}$;

\item(ii) $|{\langle}f, {f{\rangle}_2} - {\langle}G(f),
  G(f){{\rangle}_2}| {\leq}
  {\epsilon}{\|}f{\|}_{\infty}{\|}f{\|}_{\infty}$;

\item(iii)
  ${\|}{\nabla}(G(f)){\|}_2{\leq}(1+{\epsilon}){\|}{\nabla}f{\|}_2
   +{\epsilon}{\|}f{\|}_{\infty}$.

\end{description}Then if $\epsilon$ can be made arbitrarily small,
the first eigenvalue of the Jacobi operator on $M_1$ can be made
arbitrarily close to the first eigenvalue of the Jacobi operator on
${M_2}(\mu)$.

\end{lemma}

\includegraphics{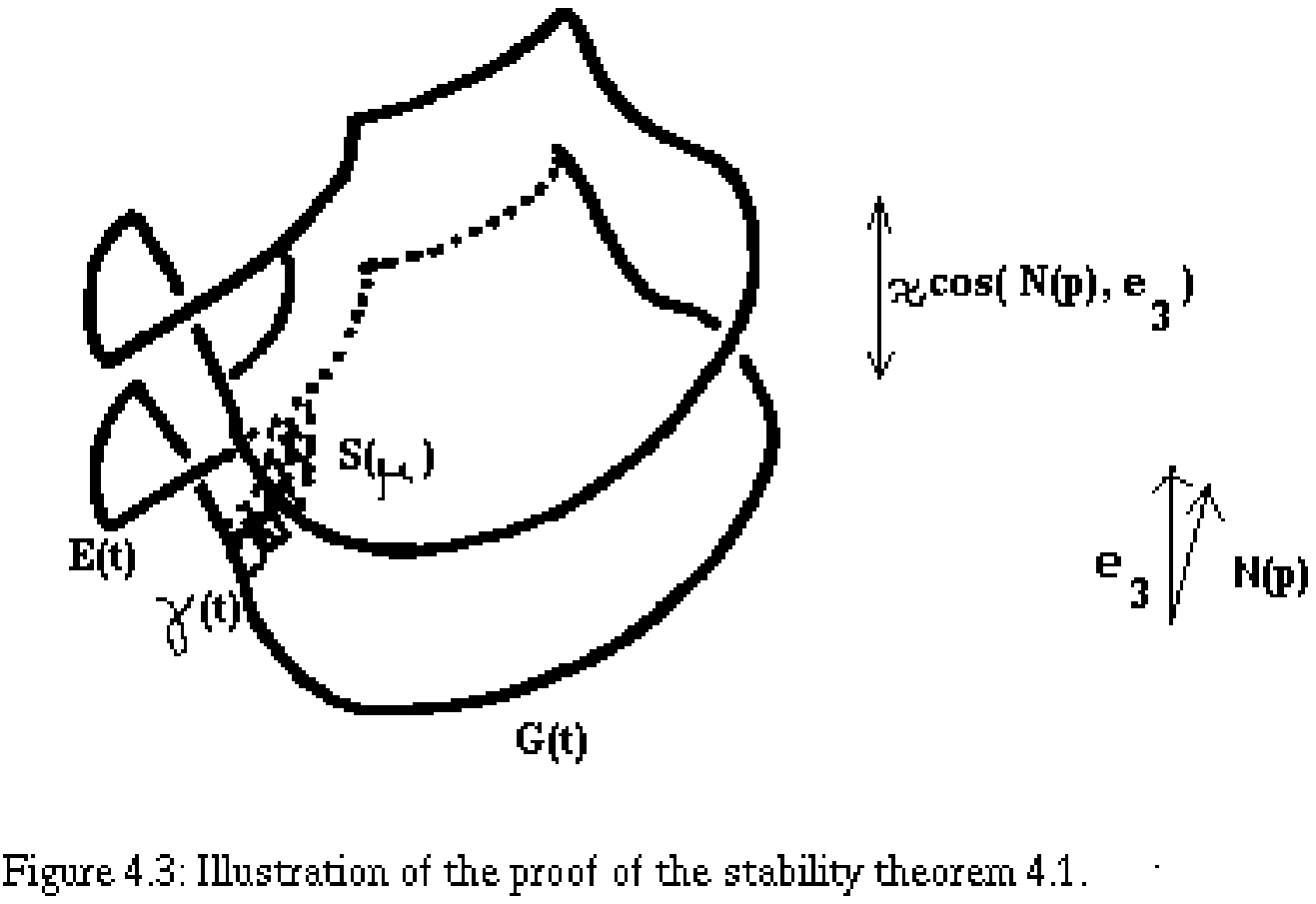}


\begin{proof}  Let us consider

  $$\lambda_1(M_1) =
  {\inf}_{\overline{f}{\in}C_0^{\infty}(M_1)}
   \frac{{\|}{\nabla}\overline{f}{\|}_2^2}{{\|}\overline{f}{\|}_2^2}.$$
  Then for each $\epsilon > 0$, if $f {\in} C_0^{\infty}(M_1)$ is a Jacobi
  vector field, one has, as observed in paragraph 1.2,
  $${\|}{\nabla}f{\|}_2^2 < (\lambda_1(M_1) +
  \epsilon){\|}f{\|}_2^2.$$ Let us choose
  ${\epsilon}=\sqrt{{Area}(S(\mu))}$.  Because of our hypotheses, $f$
  induces a $G(f) {\in} C_0^{\infty}({M_2}(\mu))$ such that condition
  (ii) above is satisfied, namely
  $$|{\,}{\|}f{\|}_2^2 - {\|}G(f){\|}_2^2{\,}| {\leq}
  \epsilon{\|}f{\|}_{\infty},$$ which implies
  $${\|}{\nabla}f{\|}_2^2 < (\lambda_1(M_1)
  +\epsilon)({\|}G(f){\|}_2+{\epsilon}{\|}f{\|}_{\infty})^2.$$
  Moreover condition (iii) above implies
  $$(\frac{{\|}{\nabla}G(f){\|}_2
    -\epsilon{\|}f{\|}_{\infty}}{1+\epsilon})^2{\leq}
  {\|}{\nabla}f{\|}_2^2,$$ which yields the conclusion:

  $$(\frac{{\|}{\nabla}G(f){\|}_2-\epsilon{\|}f{\|}_{\infty}}{1+\epsilon})^2<
  (\lambda_1(M_1)
  +\epsilon)({\|}G(f){\|}_2+{\epsilon}{\|}f{\|}_{\infty})^2.$$ If
  $\epsilon$ can be chosen arbitrarily small at the beginning, when
  $\epsilon \to 0$ one has
  $${\|}{\nabla}G(f){\|}_2^2 {\leq}\lambda_1(M_1){\|}G(f){\|}_2^2,$$
  namely $\lambda_1({M_2}(\mu)){\leq}\lambda_1(M_1)$.  Clearly in our
  case we can exchange the roles of $M_1$ and ${M_2}(\mu)$, since
  ${M_2}(\mu)$ is contained in $M_1$, and have that as ${\mu}\to 0$,
  ${\lambda}_1(M_1)\to {\lambda}_1({M_2}(\mu))$, and vice-versa.  In
  particular, we have that if $M_1$ is unstable, so is ${M_2}(\mu)$,
  and vice-versa.
\end{proof}

In the next lemma we will show that, for the surfaces $M_1$ and
${M_2}(\mu)$ defined above, the construction of $G(f)$ having the
properties required in Lemma $4.2$ is possible.

\begin{lemma}  Let $M_1 = {\Sigma}(t)$, and 
  ${M_2}(\mu) = {{\Sigma}(t)}{\backslash}S(\mu)$, ${\mu} {\approx} 0$,
  and suppose that $f {\in} C_0^{\infty}(M_1)$ defines a Jacobi vector
  field on $M_1$. Then it is possible to define a not identically zero 
  function $G(f) {\in} C_0^{\infty}({M_2}(\mu))$, in such a way that the
  conditions (i), (ii), (iii) of Lemma $3.3$ are satisfied.
\end{lemma}

\begin{proof}
  Let $f {\in} C_0^{\infty}(M_1)$ be the Jacobi vector field on $M_1$,
  and without loss of generality let us suppose that the maximum of
  $f$ is attained at some point belonging to $\gamma$ (otherwise the
  proof of the lemma is still valid, as one can easily see: this hypothesis
  takes care of the ``worst possible case'').  Moreover
  let $\phi$ be a bump function of class $C^{\infty}$ on a thin strip
  $S$ containing $\gamma$ and having area less than $\delta$ such that
  $\phi$ is constantly equal to $1$ on $M_1{\backslash}S$, constantly
  equal to zero on $S'{\subset}S$ ($S'$ is a strip containing $\gamma$
  and contained in $S$), and such that
  $|{\nabla}{\phi}| {\leq} {\frac{2}{{\delta}^{1/4}}}$.  Define now
  $G(f) = {\phi}f$.  We get:
\vskip .5cm

\begin{description}

\item{(i)} ${\|}f{\|}_{\infty} < 2{\|}G(f){\|}_{\infty}$, because 
    of the property shown in lemma $3.2$.

\item{(ii)} $|{\langle}f, f{\rangle}-{\langle}G(f), G(f){\rangle}|
  {\leq} {\delta}^2{\|}f{\|}_{\infty}{\|}f{\|}_{\infty}$, namely

\begin{eqnarray*}
  |{\|}f{\|}_2^2-{\|}{\phi}f{\|}_2^2| & = &
  |{\int}_{E{\cup}P}f^2-{\int}_{E{\cup}P}{\phi}^2{f}^2|\\ & \leq &
  {\int}_S|f^2-{\phi}^2f^2|\\ & = & {\int}_S|f^2||1-{\phi}^2|\\ &
  {\leq} &
  {\|}f^2{\|}_{\infty}{{\delta}^2}{\|}1-{\phi}^2{\|}_{\infty}\\ &
  {\leq} & {{\delta}^2}{\|}f{\|}_{\infty}^2.
\end{eqnarray*}

\item{(iii)} ${\|}{\nabla}(G(f)){\|}_2{\leq}(1+{\epsilon}){\|}{\nabla}f{\|}_2
   +{\epsilon}{\|}f{\|}_{\infty}$, since

\begin{eqnarray*} {\|}{\nabla}({\phi}f){\|}_2^2 & = &
    {\int}_{M_2}|{\nabla}({\phi}f)|^2\\ & \leq &
    {\int}_{M_2}(|{\nabla}{\phi}| |f|+|\phi| |{\nabla}f|)^2\\ & {\leq}
    & {\int}_{S}|{\nabla}{\phi}|^2 |f|^2+2{\int}_{S}|{\nabla}{\phi}|
    |f| |{\nabla}f|+{\int}_{M_2}|{\nabla}f|^2\\ & {\leq} &
    {\|}({\nabla}{\phi})^2
    (f)^2{\|}_1+2{\frac{2}{{\delta}^{\frac{1}{4}}}}{\int}_{S}|f|
    |{\nabla}f|+{\int}_{M_2}|{\nabla}f|^2\\ & {\leq} &
    {\|}({\nabla}{\phi})^2{\|}_2{\|}f^2{\|}_2
  +{\frac{4}{{\delta}^{\frac{1}{4}}}}{\|}(f)({\nabla}f){\|}_1+{\|}{\nabla}f{\|}_2^2\\ 
    & {\leq} &
    ({\int}_{S}|{\nabla}{\phi}|^4)^{\frac{1}{2}}({\int}_{S}|f|^4)^{\frac{1}{2}}
     +{\frac{4}{{\delta}^{\frac{1}{4}}}}{\|}f{\|}_2{\|}{\nabla}f{\|}_2
       +{\|}{\nabla}f{\|}_2^2\\ 
    & {\leq} &
    ({\delta}^2
     {\frac{16}{\delta}})^{\frac{1}{2}}({\delta}^2{\|}f{\|}_{\infty}^4)^{\frac{1}{2}}
     +{\frac{4}{{\delta}^{\frac{1}{4}}}}{\|}f{\|}_2{\|}{\nabla}f{\|}_2
     +{\|}{\nabla}f{\|}_2^2\\ 
    & {\leq} &
    (16{\delta})^{\frac{1}{2}}({\delta}^2{\|}f{\|}_{\infty}^4)^{\frac{1}{2}}
     +{\frac{4}{{\delta}^{\frac{1}{4}}}}({\delta}^2{\|}f{\|}_{\infty}^2)^{\frac{1}{2}}{\|}
      {\nabla}f{\|}_2+{\|}{\nabla}f{\|}_2^2\\ 
    & = &
    4{\delta}^{\frac{3}{2}}{\|}f{\|}_{\infty}^2
    +{\frac{4{\delta}}{{\delta}^{\frac{1}{4}}}}{\|}f{\|}_{\infty}{\|}{\nabla}f{\|}_2
    +{\|}{\nabla}f{\|}_2^2\\ 
    & = & (2{\delta}^{\frac{3}{4}}{\|}f{\|}_{\infty}+{\|}{\nabla}f{\|}_2)^2.
\end{eqnarray*} 
\end{description}

The assertion hence follows by choosing ${\epsilon} {\leq}
{\min}{\{}{{\delta}^2}, 2{\delta}^{\frac{3}{4}}{\}}$.
\end{proof}

In the proof of theorem $4.1$ we will denote, with the same notation
adopted previously, ${\gamma} =
{\partial}E_{t_n}{\cap}{\partial}P_{t_n}$.  We will prove the theorem
here by taking $M_1 = {\Sigma}_{t_n}$, and ${M_2}(\mu) =
{\Sigma}_{t_n}{\backslash}S(\mu)$, which we suppose to be stable.
However, since the number of crossing points is finite, the proof also
holds if we take ${M_2}({\mu}_1, \ldots , {\mu}_{n'}) =
{\Sigma}_{t_n}{\backslash}{\bigcup}_{J=1}^{n'}S({\mu}_j)$,
${n'}{\leq}n$, which we assume to be stable, and $S(\mu_j) =
{\partial}(E_{t_n}(j)){\cap}{\partial}(P_{t_n})$, where $E_{t_n}(j)$
is a sufficiently small neighborhood of the crossing point $p_j$ in
${\Sigma}_{t_n}$.

Let us now give a proof of theorem $4.1$, which we restate for easy
reference.

{\bf Theorem 4.1.} Let ${\Sigma}(t)$ be an embedded minimal surface
  with boundary ${\alpha}{\cup}{\beta}(t)$.  Suppose that
  ${\Sigma}(t)$ is described as in theorem $3.6$.  Then, for $t$
  sufficiently close to $0$, ${\Sigma}(t)$ is stable.

\begin{proof}  Suppose that the assertion stated in the theorem is false.
  Then there would exist a sequence $t_n\to 0$, corresponding to which
  there would be a sequence of unstable minimal surfaces
  ${\Sigma}(t_n)$ with boundary ${\alpha}{\cup}{\beta}(t_n)$, and
  described according to theorem $3.6$.  Hence there would exist a
  sequence of Jacobi vector fields $f_{t_n}$ defined on
  ${\Sigma}(t_n)$, all having the property proven in lemma $4.2$.
  Therefore, by induction and by lemmas $4.2$, $4.3$, $4.4$, it would be
  possible to fix a positive $\mu$ such that in a ${\mu}$-neighborhood
  (strip) of $\gamma$ one could define a bump function having bounded
  gradient which, when multiplied by $f_{t_n}$ would yield a new function
  $G(f_{t_n})$ defined on the stable part given by ${M_2}(\mu)$ (the
  fact that this would be possible for each $n$ follows from lemma
  $3.2$).  But the Rayleigh quotient associated to such a function can
  be made (dependently on $\mu$) arbitrarily close to a number which
  is strictly less than $2$, thus producing a contradiction.  For the
  sake of completeness, let us notice here that lemma $4.2$ is of
  fundamental importance in this proof, because it ensures that none
  of the functions $G(f_n)$ is the function identically equal to zero.
\end{proof}

\setcounter{chapter}{5}
\setcounter{section}{-1}
\setcounter{theorem}{0}
\setcounter{lemma}{0}
\setcounter{definition}{0}
\setcounter{claim}{0}
\setcounter{corollary}{0}
\setcounter{condition}{0}
\setcounter{question}{0}
\setcounter{example}{0}
\setcounter{remark}{0}
\section{Proof of the Main Theorem}

In this section we will put together the facts proven so far, to give
a proof of the main theorem, and to give the exact number of 
${\Sigma}(V, t)$ for $t$ sufficiently small.

\section{Uniqueness}

Let $Slab (t)$ be the slab determined by the planes $P_0$ and 
$P_t$.  Let ${\cal V}({\alpha}, {\beta})$ be the finite collection 
of varifolds determined by $\alpha{\cup}{\beta}$.  By the results 
proved in the previous sections, we know that for all values of $t$ 
sufficiently small, each varifold $V$ in ${\cal V}({\alpha}, {\beta})$ 
determines at least one compact stable embedded minimal surface 
${\Sigma}(V, t)$ in ${\cal S}(t)$ bounded by
${\alpha}{\cup}{\beta}(t)$.  We now show:

\begin{theorem}  The natural correspondence between 
  ${\cal V}({\alpha}, {\beta})$ and ${\cal S}(t)$ is a well defined
  and one-to-one correspondence, in the sense that to each varifold in
  ${\cal V}({\alpha}, {\beta})$ corresponds one and only one minimal
  surface in ${\cal S}(t)$, for $t$ sufficiently small.
\end{theorem}

\begin{proof} The theorem will be proven by contradiction.  

Let ${\{}\overline{\Sigma^1}(t) \colon {0{\leq}t{\leq}{t^{\#} {\}}}}$ and
${\{}\overline{\Sigma^2}(t) \colon {0{\leq}t{\leq}{t^{\#} {\}}}}$ be two
distinct families (for $t^{\#}$ sufficiently close to $0$) of minimal
surfaces with the same boundary ${\alpha}{\cup}{\beta}(t)$, existing
for all $0{\leq}t{\leq}{t^{\#}}$, and having the same limit varifold
$V$.  Let $M_1(t)$ be the unbounded connected component of the region
$T$ of space given by the collection of points in the slab which are
``outside'' of $\overline{\Sigma^1}(t){\cup}\overline{\Sigma^2}(t)$,
$R$ be the union of the bounded connected regions given by the points
``in between'' $\overline{\Sigma^1}(t)$ and $\overline{\Sigma^2}(t)$,
and $M_2(t)$ be the closure of $T{\backslash}R$.  Notice that $M_1(t)$
contains the truncated cylinder above $P{\backslash}V$, and that
$M_2(t)$ is contained in a small neighborhood of the truncated
cylinder above $V$.  Notice that ${\partial}({M_1}{\cup}{M_2})$
strictly contains
$\overline{\Sigma^1}(t){\cup}\overline{\Sigma^2}(t)$.  Since
$\overline{\Sigma^1}(t)$ and $\overline{\Sigma^2}(t)$ are stable by
theorem $3.1$, applying the existence theorem by Meeks and Yau stated
in section $3$ to $M_1$, we obtain the existence of a stable minimal
surface ${\Sigma_a}(t)$ ``above''
$\overline{\Sigma^1}(t){\cup}\overline{\Sigma^2}(t)$.  By the same
theorem of Meeks and Yau, applied to the region $M_2$, we obtain the
existence of a stable minimal surface ${\Sigma_b}(t)$ ``below''
$\overline{\Sigma^1}(t){\cup}\overline{\Sigma^2}(t)$.  By
construction, ${\Sigma_a}(t)$ and ${\Sigma_b}(t)$ are disjoint from
each other.  Moreover, for $t$ sufficiently close to $0$, the two
stable minimal surfaces ${\Sigma_a}(t)$ and ${\Sigma_b}(t)$ obtained
in this way can be described as stated in theorem $2.5$, namely as
approximately helicoidal around the crossing points with multiplicity
$(0, 1, 0, 1)$, and as almost horizontal graphs away from the crossing
points with multiplicity $(0, 1, 0, 1)$, by virtue of theorem $4.1$.
Hence ${\Sigma_a}(t)$ and ${\Sigma_b}(t)$ are homeomorphic to the
$\overline{\Sigma^i}(t)$, and furthermore ${\Sigma_a}(t)$ and
${\Sigma_b}(t)$ are normal graphs above $\overline{\Sigma^i}(t)$, and
hence over each other, namely there exists a function $h {\in}
C_0^{\infty}({\Sigma_b}(t))$ such that:
\protect\begin{equation}{\Sigma_a}(t)={\Sigma_b}(t)+h{\vec N}({\Sigma_b}(t)),
\protect\end{equation}that is,
for each $q{\in}{\Sigma_a}(t)$ there exists a unique
${q'}{\in}{\Sigma_b}(t)$ such that
  $$q={q'}+h(q'){\vec N}(q'),$$ where ${\vec N}(q')$ is the normal
vector to ${\Sigma_b}(t)$ in $q'$.  Hence ${\Sigma_a}(t)$ and
${\Sigma_b}(t)$ converge to the same limiting varifold, as $t{\to}0$.
It also follows from $(5.1)$ that ${\Sigma_a}(t){\cup}{\Sigma_b}(t)$
bounds a product region, say ${\cal R}(t)$.  Since ${\Sigma_a}(t)$ and
${\Sigma_b}(t)$ are normal graphs over each other, we know that the
angle between the two normal vectors to ${\Sigma_a}(t)$ and
${\Sigma_b}(t)$ at a boundary point is strictly between $0$ and $\pi$.
Then by the minimax theorems due to Pitts and Rubinstein \cite{PiRu},
and generalized by Jost \cite{Jo} to the case of nonempty boundary,
${\alpha}{\cup}{\beta}(t)$ is also the boundary of an unstable
embedded minimal surface ${\Sigma^*}(t)$, contained in $\cal R$, for
all $t$ sufficiently small.  We now wish to show that ${\Sigma^*}(t)$
is homeomorphic to ${\Sigma_a}(t)$ and ${\Sigma_b}(t)$, and that
${\Sigma^*}(t)$ actually has the same geometric description as
${\Sigma_a}(t)$ and ${\Sigma_b}(t)$, which will imply, for example,
that $\lim_{t{\to}0}{\Sigma^*}(t)
=\lim_{t{\to}0}{\Sigma_a}(t)=\lim_{t{\to}0}{\Sigma_b}(t)=V$.  First
notice that the connected components of the complement in
${\Sigma^*}(t)$ of the union of small neighborhoods of the crossing
points, is the union of almost horizontal graphs, just like
${\Sigma_a}(t)$ and ${\Sigma_b}(t)$.  In fact, the area of each such
connected component is almost equal to the area of the corresponding
components of ${\Sigma_a}(t)$ and ${\Sigma_b}(t)$.  The proof of this
fact follows from comparing the area of ${\Sigma^*}(t)$ with the
maximum area of the family ${\Sigma_a}(t)={\Sigma_b}(t)+th{\vec
N}({\Sigma_b}(t))$, indexed by $t{\in}[0, 1]$.  The definition of
minimax implies that the area of ${\Sigma^*}(t)$ cannot be larger than
the maximum area of the family; but since every surface in the family
is a graph over $\Sigma_b(t)$, we can estimate the area of
${\Sigma^*}(t)$ by that of ${\Sigma_a}(t)$ and ${\Sigma_b}(t)$, and
the proof of the claim follows as in part $1$ of the proof of theorem
$3.6$, since the area estimates just proven allow us to apply
R. Schoen's curvature estimates.  So we know that ${\Sigma^*}(t)$
given by very flat graphs away from the crossing points.  In the
neighborhood of a crossing point with multiplicity $(0, 1, 0, 1)$,
consider the part of $\Sigma^*(t)$ bounded by four geodesic arcs,
constructed as in part $3$ of the proof of theorem $3.6$.  Such
quadrilateral region has sum of the external angles strictly less than
$4\pi$.  By Nitsche's $4{\pi}$-theorem stated in section $2$, applied
to an analytic smoothing having total curvature less than $4\pi$ of
the quadrilateral defined by the geodesic arcs, we know that this
quadrilateral region must bound a {\em stable} minimal surface which
is topologically a disk.  This means that $\Sigma^*(t)$ can be
described as in theorem $3.6$, for all $t$ sufficiently small, and
that ${\Sigma}^*(t)$ is a graph over ${\Sigma}_b(t)$.  But then, for
$t$ sufficiently close to $0$, $\Sigma^*(t)$ must be stable, by
theorem $4.1$.  This produces a contradiction and finishes the proof
of the uniqueness theorem.
\end{proof}

\includegraphics{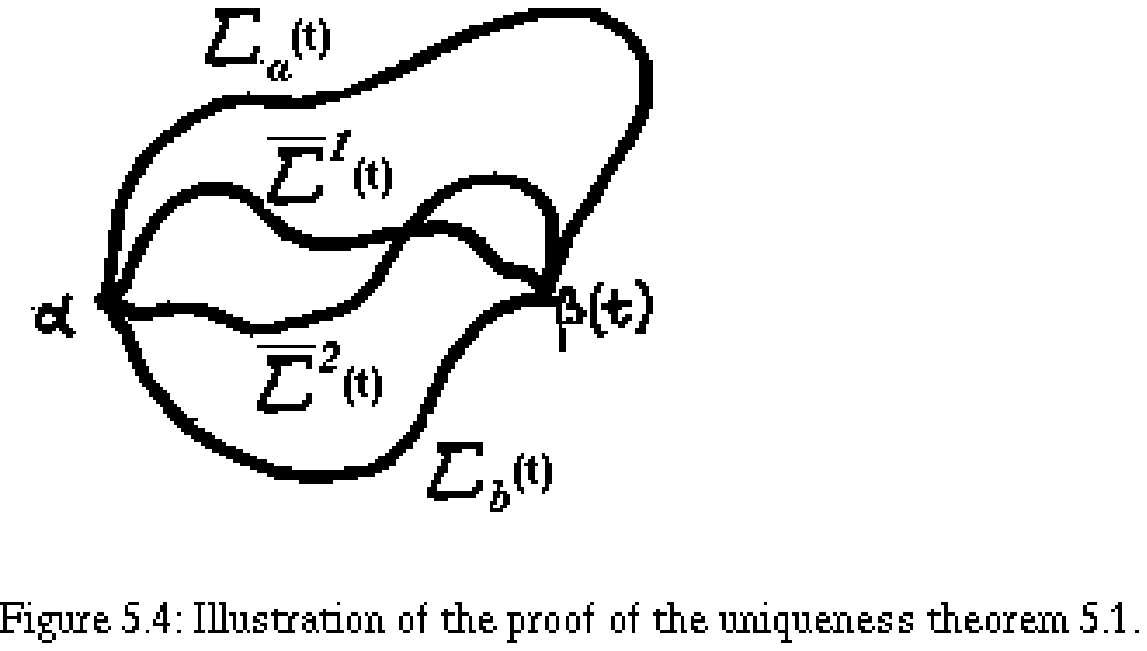}

\section{A bound on the number of ${\Sigma}(V, t)$'s}


In this section we will derive a formula which gives an upper bound for
the number of compact stable minimal surfaces bounded by a finite number
of Jordan curves in close planes of ${\R}^3$.  

By theorem $5.1$, we know that once the multiplicity of the limiting
varifold $V$ is fixed, then there is a unique stable compact minimal
surface bounded by ${\alpha}_t{\cup}\beta$, for $t$ sufficiently close
to $0$.  So the number of stable compact embedded minimal surfaces
will be determined once we are able to express exactly how many are
the possible multiplicities for the limiting varifold $V$.  In order
to understand this, let us assign to each connected component of
${P_0}{\backslash}({\alpha}\cup{\beta})$ the sign ``$+$'' or ``$-$''
in such a way that the unbounded component $C_{u}$ is given the
``$-$'' sign and adjacent components have opposite signs.  Moreover
notice that ${P_0}{\backslash}({\alpha}\cup{\beta})$ determines a
finite number of varifolds, by assigning to each of its connected
components one of the numbers $0$, $1$ or $2$ (note that the connected
components having the ``$+$'' sign can only be assigned multiplicity
one).

So the multiplicity is totally determined for the components of $V$
with the ``$+$'' sign.  The only places where there are different
multiplicities allowed are the components of $V$ with the ``$-$''
sign.  Of course we only have two choices for the multiplicity of such
components: $zero$ or $two$.  However we are not free to assign
multiplicities arbitrarily, as $(1, 2, 1, 2)$ crossing points must be
avoided.  Let $f^i_{-}$ be the number of components of $V$ having
``$-$'' sign, and contained inside ${\cal R}(\alpha)\cap{\cal
R}(\beta)$, and let $f^o_{-}$ be the number of components of $V$
having ``$-$'' sign, and contained outside ${\cal R}(\alpha)\cup{\cal
R}(\beta)$.  Then the above observations translate in the following

\begin{corollary} Once the limiting cycle $Z$, is given, the number of
  stable compact minimal surfaces ${\Sigma}(t)$ such that
  ${\partial}({\Sigma}(t)) {\to} Z$ as $t{\to}0$, is bounded above by

  $$2^{f^i_{-}} + 2^{f^o_{-}}.$$ 

\end{corollary}

\begin{remark} It would be interesting to get similar bounds 
  on the number of unstable embedded minimal surfaces bounded
  by Jordan curves in close planes.
\end{remark}


\setcounter{chapter}{6}
\setcounter{section}{-1}
\setcounter{theorem}{0}
\setcounter{lemma}{0}
\setcounter{definition}{0}
\setcounter{claim}{0}
\setcounter{corollary}{0}
\setcounter{condition}{0}
\setcounter{question}{0}
\setcounter{example}{0}
\setcounter{remark}{0}

\section{A nonexistence result}

In this section we observe that our main theorem provides some
evidence in support of a conjecture made by W. Meeks in \cite{Me1}.
The conjecture states that there are no minimal surfaces of positive
genus bounded by two convex curves in parallel planes of ${\R}^3$.
One consequence of our main theorem is:

\begin{corollary}  There exist no stable minimal surfaces of positive 
genus bounded by two convex curves in parallel planes of ${\R}^3$,
when the distance between the two planes is sufficiently small (or 
if the two planes are not parallel, but they intersect at a sufficiently
small angle).  \end{corollary}

Previous results by R. Schoen \cite{Sc1}, W. Meeks and B. White 
\cite{MeWh1} \cite{MeWh2} \cite{MeWh3}, and earlier M. Shiffman
\cite{Sh}, also supported evidence towards Meeks's conjecture.

\begin{remark}  If one could show that as $t$, the distance
between the planes, increases, the number of {\em stable} minimal
surfaces bounded by two convex curves does not increase, Meeks's
conjecture would follow, at least in the stable case.
\end{remark}

}

\end{document}